\documentclass[a4paper]{amsart}
\usepackage{sabbah_Fourier-irred}

\begin{document}
\title[Fourier-Laplace transform]{Fourier-Laplace transform of\\
irreducible regular differential systems on\\
the Riemann sphere}

\author[C.~Sabbah]{Claude Sabbah}
\address{UMR 7640 du CNRS\\
Centre de Math\'ematiques Laurent Schwartz\\
\'Ecole polytechnique\\
F--91128 Palaiseau cedex\\
France}
\email{sabbah@math.polytechnique.fr}
\urladdr{http://www.math.polytechnique.fr/cmat/sabbah/sabbah.html}
\dedicatory{\`A la m\'emoire d'Andre\"i Bolibroukh}

\begin{abstract}
We show that the Fourier-Laplace transform of an irreducible regular differential system on the Riemann sphere underlies, when one only considers the part at finite distance, a polarizable regular twistor $\mathcal{D}$-module. The associated holomorphic bundle out of the origin is therefore equipped with a natural harmonic metric with a tame behaviour near the origin.
\end{abstract}

\subjclass{Primary 32S40; Secondary 14C30, 34Mxx}

\keywords{Flat bundle, harmonic metric, twistor $\cD$-module, Fourier-Laplace transform}




\maketitle
\tableofcontents

\section*{Introduction}
One of the important results of A.~Bolibrukh was to give a positive answer to the Riemann-Hilbert problem for semisimple\footnote{that is, direct sum of irreducible objects}\addtocounter{footnote}{-1} linear representations of the fundamental group of the complement of a finite set of points on the Riemann sphere. Similarly, he proved that the Birkhoff problem has a positive answer when the corresponding system of meromorphic linear differential equations is semisimple\footnotemark, generalizing previous results of W.~Balser.

Both results can be stated in a similar way, by using the language of meromorphic bundles on the Riemann sphere: let $P=\{p_1,\dots,p_r,p_{r+1}=\infty\}$ be a non empty finite set of points on $\PP^1$ and let $M$ be a free $\cO_{\PP^1}(*P)$-module of finite rank, equipped with a connection $\nabla$; assume that there exists a basis of $M$ in which, for $i=1,\dots,r$, the connection $\nabla$ has Poincar\'e rank $m_i\geq0$ at $p_i$ (\ie the order of the pole of the matrix of $\nabla$ at $p_i$ is $m_i+1$); assume moreover that $\nabla$ has a regular singularity at $\infty$ (\ie the coefficients in the given basis of horizontal sections are multivalued holomorphic functions which grow at most polynomially at infinity); then, if $(M,\nabla)$ is \emph{irreducible} (or semisimple), there exists a basis of $M$ in which the Poincar\'e rank of $\nabla$ at $p_i$ is $m_i$ ($i=1,\dots,r$) and $\nabla$ has at most a logarithmic pole at $\infty$ (\ie $\nabla$ has Poincar\'e rank~$0$ at infinity).

After an easy preliminary reduction, the case where $m_i=0$ for all $i$ corresponds to the Riemann-Hilbert problem and the case where $r=1$ to the Birkhoff problem.

Starting from a system of linear meromorphic differential equations having only regular singularities --- or more precisely from a regular holonomic $\cD$-module --- on the Riemann sphere, one obtains, by Fourier-Laplace transform, a new system, and the Riemann-Hilbert problem for the original system transforms to the Birkhoff problem for the new one. These problems (for a given system and its Fourier-Laplace transform) are not directly related one to the other\footnote{See, however, \cite[\S V.2.c]{Bibi00}.} but, as one of both systems is semisimple if and only if the other one is so, there is a common condition under which both problems simultaneously have a positive answer.

Semisimple linear representations as above share another remarkable property: on the associated flat bundle there exists a \emph{tame harmonic metric} (\cf \cite{Simpson90} or see below). This can be expressed by using the language of \emph{polarized twistor $\cD$-modules} introduced in \cite{Bibi01}, extending a notion due to C.~Simpson \cite{Simpson97}: to such a representation is associated a unique --- up to isomorphism --- regular holonomic $\cD$-module on the Riemann sphere, which has neither submodule nor quotient module supported at a point (\cf Lemma \ref{lem:RHC} below); the previous property can be expressed as follows: this $\cD$-module underlies a polarizable regular twistor $\cD$-module, as defined in \loccit (see below).

In this article, we analyze the behaviour of polarized regular twistor $\cD$-modules on the Riemann sphere by Fourier-Laplace transform. However, we do not give any information on the behaviour at infinity in the Fourier plane, where an irregular singularity occurs.

It should be mentioned that, using different techniques, S.~Szabo \cite{Szabo04} has established a perfect Fourier correspondence in a more general situation, where an irregular singularity is allowed at infinity, but with some other more restrictive assumptions.

We refer to \cite{Bibi01} for the various results we use concerning polarizable twistor $\mathcal D$-modules.

\section{Statement of the results}\label{sec:statement}
We view the projective line $X=\PP^1$ as the union of the two affine charts $\spec\CC[t]$ and $\spec\CC[t']$ with $t'=1/t$ on the intersection, and we define $\infty$ as the point where $t'=0$. As above, let $P=\{p_1,\dots,p_r,p_{r+1}=\infty\}$ be a finite set of $r+1$ distinct points in $\PP^1$. Put $X^*=\PP^1\moins P$.

Let $(H,D_V)$ be a $C^\infty$ vector bundle with a flat connection on $X^{*\an}$. It is holomorphic with respect to the $(0,1)$-part $D''_V$ and we put $(V,\nabla)=(\ker D''_V,D'_V)$. The associated local system is $\cL\defin\ker\big[\nabla:V\to V\otimes_{\cO_{X^{*\an}}}\Omega^1_{X^{*\an}}\big]$. We denote by
$T_j$ the local monodromy of this local system at each point $p_j$ of $P$, that we view as an automorphism of the generic fibre $L$ of $\cL$ (or of $H$).

\subsection{Fourier-Laplace transform of the flat bundle}
The notion of Fourier-Laplace transform is \emph{a priori} neither defined for holomorphic bundles with connection on~$X^*$ nor for holomorphic vector bundles on $\PP^1$ equipped with a meromorphic connection. It is defined for algebraic $\cD$-modules on the affine line $\Afu=\PP^1\moins\{\infty\}$.

We denote by $\Clt$ the Weyl algebra in dimension $1$, that is, the quotient of the free algebra generated by $\CC[t]$ and $\CC[\partial_t]$, by the relation $[\partial_t,t]=1$ (see \eg \cite{Borelal87,Bibi90}). Let $M$ be a holonomic $\Clt$-module. We say that $M$ is a \emph{minimal extension} if it has neither submodules nor quotient modules supported on some point in $\Afu$. The following is well-known.

\begin{lemme}[Riemann-Hilbert correspondence]\label{lem:RHC}
The functor which associates to any holonomic $\Clt$-module having singularities at most at $P$ its restriction to $X^{*\an}$ induces an equivalence between the category of holonomic $\Clt$-modules which have regular singularities (included at infinity) and are a minimal extension (morphisms being all morphisms of $\Clt$-modules) and the category of flat holomorphic bundles on $X^{*\an}$.\qed
\end{lemme}

Given a holomorphic bundle with connection $(V,\nabla)$ on $X^{*\an}$, we denote by $M$ the regular holonomic $\Clt$-module associated to $(V,\nabla)$ by Lemma \ref{lem:RHC}. The Fourier-Laplace transform $\wh M$ of $M$ is the $\CC$-vector space $M$ equipped with an action of the Weyl algebra $\Cltau$ in the variable $\tau$ defined by the formula
\begin{equation}
\tau\cdot m=\partial_t m,\quad \partial_\tau m=-tm.
\end{equation}
We will denote by $\wh X=\Afuh$ the affine line with coordinate $\tau$ and by $\wh\PP^1$ the corresponding projective line. It is known that $\wh M$ has a single singularity at finite distance, namely at $\tau=0$, which is regular. However, it has in general an irregular singularity at $\tau=\infty$ (see \eg \cite{Malgrange91} for general results on Fourier transform of holonomic $\Clt$-modules). We therefore get a holomorphic bundle with connection $(\wh V,\wh\nabla)$ on $\wh X^*=\Afuh\moins\{0\}=\wh\PP^1\moins\wh P$, where $\wh P=\{0,\infty\}$. The associated flat $C^\infty$ bundle is denoted by $(\wh H,D_{\wh V})$. However, $\wh M$, having an irregular singularity at infinity, cannot be recovered from the only datum of $(\wh V,\wh\nabla)$, in general.

Denote by $d$ the rank of $H$. Let us indicate how to compute the rank $\wh d$ of $\wh H$. It is known that, for any $\tau_o\neq0$, the $\CC$-linear morphism
\[
M\To{\partial_t-\tau_o}M,
\]
is injective and that its cokernel is a finite dimensional vector space: this is the fibre of $\wh V$ at the point $\tau_o$. Its dimension, which is the total number $\mu$ of ``vanishing cycles'' of $M$ at $P\moins\{\infty\}$ (\ie the sum of multiplicities of the characteristic variety of $M$ along its components $T^*_{p_i}X$, see \eg \cite[\S4, p.~66]{Malgrange91}), is easily computed here by the formula  (\cf \cite[Prop.~(1.5), p.~79]{Malgrange91})
\begin{equation}\label{eq:dhat}
\wh d=\mu=r\cdot d-d_1, \quad r=\card P -1,\quad d_1\defin\sum_{j=1}^r\dim\ker(T_j-\id).
\end{equation}
More precisely, $\CC[\tau,\tau^{-1}]\otimes_{\CC[\tau]}\wh M$ is a free $\CC[\tau,\tau^{-1}]$-module of rank $\mu$.

\subsection{Fourier-Laplace transform of an irreducible bundle with connection}

Assume now that $(H,D_V)$, is irreducible, or, equivalently, that there is no nontrivial subspace of $L$ invariant under all $T_j$ ($j=1,\dots,r+1$). This is equivalent to saying that the associated $M$ is irreducible as a $\Clt$-module (easy). This, in turn, is equivalent to the irreducibility of $\wh M$ (easier). However, $\wh M$ being irregular at infinity, this does not imply the irreducibility of $(\wh H,D_{\wh V})$ (see below). The rank~$\wh d$ of~$\wh H$ is positive unless $(V,\nabla)=(\cO_X,d)$, as we assume irreducibility. In the following, we implicitly assume that $\wh d>0$.

The bundle $(\wh H,D_{\wh V})$ is determined by its monodromy $\wh T_0$ around $\tau=0$. The Jordan structure of $\wh T_0$ is determined by that of the monodromy $T_\infty$ of $(H,D_V)$ around $\infty$ (we denote below by $A_\lambda$ the restriction of the endomorphism $A$ to the generalized eigenspace corresponding to the eigenvalue $\lambda$): this is the content of the following lemma, which can be deduced from Proposition~8.4.20 in \cite{Bibi01}.

\begin{lemme}\label{lem:monoinfty}
Assume that $(H,D_V)$ is irreducible. Then,
\begin{enumerate}
\item
for any $\lambda\neq1$, the Jordan structures of $\wh T_{0,\lambda}$ and $T_{\infty,\lambda}$ coincide;
\item
any Jordan block of size $k\geq1$ of $T_{\infty,1}$ gives rise to a Jordan block of size $k+1$ of $\wh T_{0,1}$.
\item
The remaining Jordan blocks of $\wh T_{0,1}$ have size one.\qed
\end{enumerate}
\end{lemme}

The semisimplicity of $(\wh H,D_{\wh V})$ (decomposability as the direct sum of irreducible bundles with connection) is equivalent to that of $\wh T_0$. This occurs (still assuming that $(H,D_V)$ is irreducible) if and only if $T_\infty$ is semisimple and $1$ is not an eigenvalue of~$T_\infty$.

\subsection{Fourier-Laplace transform of a bundle with connection and a Hermitian metric}\label{subsec:FLflatherm}

We assume first that the monodromy representation associated with the local system $\cL=\ker\nabla$ is \emph{unitary}. In other words, we assume that there exists a \emph{$D_V$-flat Hermitian metric} $h$ on $H$. In particular,
\begin{enumerate}
\item
the local system $\cL$ is an orthogonal direct sum of irreducible local systems; in the following, we will assume that it is \emph{irreducible}; the case of the constant local system of rank one being trivial, we will also assume that the local system $\cL$ is not constant.
\item
the local monodromy $T_j$ of the local system at each point $p_j$ of $P$, being unitary, is semi-simple, its eigenvalues are roots of unity.
\end{enumerate}

We may ask whether $\wh H$ comes equipped naturally with a Hermitian metric and if this metric is $D_{\wh V}$-flat or not. Notice that flatness would imply semisimplicity and unitarity of the monodromy $\wh T_0$. Under the previous assumption, this holds if and only of $1$ is not an eigenvalue of $T_\infty$, by Lemma \ref{lem:monoinfty}.

If we only assume that $(H,D_V)$ is irreducible but not necessarily unitary, there exists however, by \cite{Simpson90}, a unique \emph{tame harmonic Hermitian metric} $h$ on $(H,D_V)$, which is therefore the natural metric to consider. Such a metric also exists if $(H,D_V)$ semisimple, but may not be unique.

\subsubsection*{Tame harmonic metric}
Let us recall the definition of such a metric. We can fix the choice of a metric connection on $H$, that we denote $D_E$, by the following constraint: if we put $\theta_E=D_V-D_E$, that we decompose into its $(1,0)$ and $(0,1)$ parts as $\theta_E=\theta'_E+\theta''_E$, then we impose that the $h$-adjoint of $\theta'_E$ is $\theta''_E$. As $X$ has dimension one, the bundle $E=\ker D''_E$ is holomorphic on $X^{*\an}$ and $\theta'_E$ satisfies the \emph{Higgs condition} $\theta'_E\wedge\theta'_E=0$. 

The triple $(H,D_V,h)$ (with $D_V$ flat) is said to be \emph{harmonic} if the Higgs field is holomorphic on $E$, that is, if $\theta'_E:E\to E\otimes_{\cO_{X^{*\an}}}\Omega^1_{\cO_{X^{*\an}}}$ is holomorphic.

Following \cite{Simpson90}, we say that such a triple $(H,D_V,h)$ is \emph{tame} if the eigenvalues of the Higgs field (which are multivalued holomorphic one-forms), have pole of order at most one at each point of $P$.

One may ask whether 
$(\wh H,D_{\wh V})$ also carries such a harmonic metric, with a tame behaviour at $\tau=0$. We give a positive answer in Corollary \ref{cor:Fmain}.

\subsubsection*{Twistor $\cD$-modules}
We will use the language of twistor $\cD$-modules of \cite{Bibi01}, which we refer to. Let us quickly recall some basic definitions.

We still denote by $X$ the Riemann sphere and we put $\cX=X\times\CC$, using the coordinate $\hb$ on the factor $\CC$. We will also denote by $\bS$ the circle $\{\module{\hb}=1\}$. We denote by $\cD_X$ the sheaf of holomorphic differential operators on $X$ and we consider the sheaf $\cR_\cX$ on $\cX$ of $\hb$-differential operators: in any local coordinate $x$ on $X$, $\cR_{\cX}=\cO_{\cX}\langle\partiall_x\rangle$, with $\partiall_x=\hb\partial_x$. In particular, we have $\cD_X=\cR_\cX/(\hb-1)\cR_\cX$. 

We consider the category $\RTriples(X)$ having objects of the form $(\cM',\cM'',C)$, where $\cM',\cM''$ are coherent $\cR_\cX$-modules and $C$ is a sesquilinear pairing between them: for any $\hb_o$ in $\bS$, it is a pairing between $\cM'_{\hb_o}$ and the conjugate of $\cM''_{-\hb_o}$ taking values in the sheaf $\Db_{X_\RR}$ of distributions on $X$; it has to be linear with respect to the action of holomorphic differential operators on $\cM'_{\hb_o}$ and that of antiholomorphic differential operators on the conjugate of $\cM''_{-\hb_o}$, both kinds of differential operators acting naturally on distributions on~$X$; last, this pairing has to be continuous with respect to $\hb\in\bS$. We will view $C$ as a sesquilinear pairing $\cMS'\otimes_{\cO_{\cX|\bS}}\ov{\cMS''}\to\Dbh{X}$, where $\Dbh{X}$ denotes the sheaf on $X\times\bS$ of distributions which are continuous with respect to $\hb\in\bS$ and where the conjugation is opposite to the usual one on $\bS$ (see \cite[\S1.5.a]{Bibi01}).

In \cite{Bibi01} is defined a notion polarized regular twistor $\cD$-module of weight $w\in\ZZ$ on~$X$: such an object takes the form $(\cT,\cS)$, where $\cT=(\cM',\cM'',C)$ is as above and $\cS$ (called the \emph{polarization}) consists of two isomorphisms $\cM''\isom\cM'$. Some axioms are required, that we do not recall here (\cf \loccit, Chap.~4). In most arguments, one is able to reduce to the case where the weight $w$ is $0$, $\cM'=\cM''$ and both isomorphisms in $\cS$ are equal to $\id$. We denote then such an object simply by $(\cM,\cM,C,\id)$ or even by $(\cM,\cM,C)$.

Given such an object of weight $0$, we obtain a $\cD_X$-module by considering \hbox{$\cM/(\hb-1)\cM$}. Moreover, by restricting to $X^{*\an}$, we obtain a holomorphic bundle $V$ with connection $\nabla$. Last, the axioms imply that the sesquilinear pairing $C$ allows one to define a metric~$h$ on the $C^\infty$ bundle $H$ associated with $V$ and that this metric is \emph{harmonic}. The regularity assumption implies moreover that this metric is \emph{tame}. More precisely, using results of \cite{Simpson90} and \cite{Biquard97}, we prove in \cite[Chap.~5]{Bibi01} that the category of polarized regular twistor $\cD_X$-modules is equivalent to that of tame harmonic bundles $(H,D_V,h)$ on $X^{*\an}$ with some parabolic structure (called ``of Deligne type'' in \loccit). The latter category is shown to be equivalent to that of semisimple bundles with flat connection $(H,D_V)$ on $X^{*\an}$ by C.~Simpson in \cite{Simpson90}.

The Fourier-Laplace transform of an object $(\cM',\cM'',C,\cS)$ is defined in \cite[Chapter~8]{Bibi01} (see also below). It is a triple $(\wh\cM',\wh\cM'',\wh C,\wh\cS)$ on $\Afuhan\times\CC$. The main result of this article is:

\begin{theoreme}\label{th:Fmain}
If $(\cM,\cM,C,\id)$ is a polarized regular twistor $\cD_X$-module of weight~$0$, then $(\wh\cM,\wh\cM,\wh C,\id)$ a polarized regular twistor $\cD$-module of weight~$0$ on~$\Afuhan$.
\end{theoreme}

It is straightforward to extend this statement to polarized regular twistor $\cD$-module of weight $w$. Notice that part of the theorem has yet been proved in Theorem 8.4.1 of \loccit, namely, the condition on nearby cycles at $\tau=0$. We will therefore be mainly interested in the behaviour at $\tau_o$ with $\tau_o\neq0$. The ``fibre'' at $\tau=\tau_o\neq0$ of $(\wh\cM,\wh\cM,\wh C)$ is obtained from that at $\tau=1$ by a rescaling, that is, by a preliminary change of variable $t\mto t/\tau_o$, as the kernel of the Fourier-Laplace transform is $e^{-t\tau/\hb}$.

As a consequence, we get:

\begin{corollaire}\label{cor:Fmain}
If the flat bundle $(H,D_V)$ is semisimple, the Fourier-Laplace transform $(\wh H,D_{\wh V})$ on $\wh X^*$ carries a harmonic metric with a tame behaviour at $\tau=0$. 
\end{corollaire}

\section{Exponential twist of harmonic bundles and twistor $\cD$-modules}
Let us recall a basic correspondence given in \cite[\S8.1.b]{Bibi01}. We keep notation of \S\ref{sec:statement} but we now fix $\tau_o=1$.

\subsection{Exponential twist of smooth twistor structures}
We start with $(H,h,D_V)$ on $X^*$ and we rescale the metric $h$ and twist the connection $D_V$, defining thus
\begin{align*}
\Fou D_V&=e^t\circ D_V\circ e^{-t},\quad \text{\ie } \Fou D'_V=D'_V-dt,\ \Fou D''_V=d'',\\
\Fou h&=e^{2\reel t}h.
\end{align*}

Recall that, using definitions in \cite{Simpson90,Simpson92}, if the triple $(H,D_V,h)$ is harmonic on $X^*$, then so is the triple $(H,\Fou D_V,\Fou h)$. The Higgs field is given by the formulas
\[
\Fou\theta'_E=\theta'_E-dt,\quad\Fou\theta''_E=\theta''_E-d\ov t,
\]
and the $\Fou h$-metric connection $\Fou D_E=\Fou D'_E+\Fou D''_E$ by
\[
\Fou D_E=e^{-\ov t}\circ D_E\circ e^{\ov t},\quad\text{\ie } \Fou D'_E=D'_E,\ \Fou D''_E= D''_E+ d\ov t.
\]

The exponential twist exists at the level of smooth twistor structures. As in \cite{Bibi01}, we denote by $\cC^{\infty,\an}_{\cX^*}$ the sheaf on $\cX^*$ of $C^\infty$ functions which are holomorphic with respect to $\hb$. Consider the $\cC^{\infty,\an}_{\cX^*}$-module\footnote{In \cite{Bibi01}, this module is denoted by $\cH$; here, we stress upon the analytic dependence on $\hb$.} $\cH^{\an}=\cC^{\infty,\an}_{\cX^*}\otimes_{\pi^{-1}\cC^\infty_{X^*}}\pi^{-1}H$, equipped with the $d''$ operator
\begin{equation}
\Fou \ccD''_{\hb}=\Fou D''_E+\hb\Fou\theta''_E=\ccD''_{\hb}+(1-\hb)d\ov t.
\end{equation}
We get a holomorphic subbundle $\FcH'=\ker \Fou\ccD''_{\hb}\subset\cH^\an$ equipped with a $\hb$-connection $\Fou\ccD'_{\hb}=\hb \Fou D'_E+\Fou\theta'_E=\ccD'_\hb-dt$. We put $$\Fou\ccD_\hb=\Fou\ccD'_{\hb}+\Fou\ccD''_{\hb}=\ccD_\hb-dt+(1-\hb)d\ov t.$$
Moreover, if $\pi:\cX^*=X^*\times\CC\to X^*$ denotes the natural projection, $\cH^\an$ can be equipped with the metric $\pi^*h$ or the metric $\pi^*\Fou h$. These metrics are constant with respect to $\hb$. We will also consider the metric $e^{2\reel(\hb\ov t)}\pi^*h$, which varies with $\hb$.

We have an isomorphism of locally free $\cC^{\infty,\an}_{\cX^*}$-modules with metric and $\hb$-conn\-ec\-tion:
\begin{equation}\label{eq:isoFF}
(\cH^\an,\pi^*\Fou h,\Fou\ccD_\hb)\To{{}\cdot e^{(1-\hb)\ov t}}(\cH^\an,e^{2\reel(\hb\ov t)}\pi ^*h,\ccD_\hb-dt).
\end{equation}
This isomorphism sends the holomorphic subbundle $\Fou\cH'$ to $\cH'=\ker\ccD_\hb''$.

It will also be useful to have a model where the metric is $\pi^*h$. This model is defined on the sheaf $\cC^\infty_{\cX^*}$ and not on $\cC^{\infty,\an}_{\cX^*}$. Put $\cH=\cC^{\infty}_{\cX^*}\otimes_{\pi^{-1}\cC^\infty_{X^*}}\pi^{-1}H$. We have an isomorphism
\begin{equation}\label{eq:isoFFR}
(\cH,\pi^*\Fou h,\Fou\ccD_\hb)\To{{}\cdot e^{\ov t-2i\im(\hb\ov t)}}(\cH,\pi ^*h,\ccD_\hb-(1+\module{\hb}^2)dt).
\end{equation}
This isomorphism is not defined over $\cC^{\infty,\an}_{\cX^*}$.

\subsection{Exponential twist in $\RTriples(X^*)$} Recall the following definitions (\cf \cite[\S8.1.a]{Bibi01}). Let $\cM$ be a left $\cR_\cX$-module, \ie a $\cO_\cX$-module with a flat relative meromorphic connection $\nabla_{\cX/\CC}$ (relative is with respect to $\hb$, that is, there is no differentiation with respect to $\hb$). Denote by $\cM_\loc$ the localized module along $P$, that is, $\cM_\loc=\cO_\cX[*(P\times\CC)]\otimes_{\cO_\cX}\cM$. The twisted $\cR_\cX$-module $\FcM_\loc=\cM_\loc\otimes \ccE^{-t/\hb}$ is defined as the $\cO_\cX$-module $\cM_\loc$ equipped with the twisted connection $e^{t/\hb}\circ\nabla_{\cX/\CC}\circ e^{-t/\hb}$.

Let $C:\cMS'\otimes_{\cO_{\cX|\bS}}\ov{\cMS''}\to\Dbh{X}$ be a sesquilinear pairing. If the restriction of $(\cM',\cM'',C)$ to $\cX^*$ is a smooth twistor structure, the restriction to $X^*\times\bS$ of $C$ takes values in $\cC^{\infty,\an}_{\cX^*}$ and the extension $C_\loc$ of $C$ to $\cM_{\loc|\bS}'\otimes_{\cO_{\cX|\bS}}\ov{\cM_{\loc|\bS}''}$ takes values in the extension of $\cC^{\infty,\an}_{\cX^*}$ made of functions on $\cX^*$ which can be extended as distributions which are continuous relatively to $\hb\in\bS$. Moreover, if we assume that $(\cM',\cM'',C)$ underlies a polarized regular twistor $\cD$-module, then, using \cite[(5.3.3)]{Bibi01}, we obtain that $C_\loc$ takes values in the extension of $\cC^{\infty,\an}_{\cX^*}$ made of functions on $\cX^*$ which have moderate growth near each puncture~$P$, locally uniformly with respect to $\hb\in\bS$.

Notice that, for $\hb\in\bS$, \hbox{$\hb\ov t-t/\hb$} is purely imaginary. Then, under the previous assumption, $\Fou C_\loc\defin \exp(\hb\ov t-t/\hb)C_\loc$ is a sesquilinear pairing on $\FcM_{\loc|\bS}'\otimes_{\cO_{\cX|\bS}}\ov{\FcM_{\loc|\bS}''}$, taking values in the same sheaf of functions with moderate growth.

To a harmonic bundle $(H,h,D_V)$ on $X^*$ is attached the smooth twistor structure $(\cH',\cH',\pi^*h_{\cH'_{|\bS}\otimes\ov{\cH'_{|\bS}}})$, where $\cH'\subset\cH^\an$ is $\ker\ccD''_\hb$ equipped with the $\cR_{\cX^*}$-structure given by the $\hb$-connection $\ccD'_\hb$.

It can be exponentially twisted as an object of $\RTriples(X^*)$: the result is $(\cH',\cH',\exp(\hb\ov t-t/\hb)\pi^*h_{\cH'_{|\bS}\otimes\ov{\cH'_{|\bS}}})$, where $\cH'$ is equipped with the $\cR_{\cX^*}$-structure given by the $\hb$-connection $\ccD'_\hb-dt$.

The isomorphism \eqref{eq:isoFF} identifies it with the smooth twistor structure attached to $(H,\Fou h,\Fou D_V)$ (\cf \cite[Lemma~8.1.2]{Bibi01}).

\subsection{Exponential twist in $\RTriples(X)$}
Let $\cM$ be a left $\cR_\cX$-module. Denote by $\wt\cM$ its localization \emph{at infinity only}. Then $\FcM$ is defined as the twisted $\cR_\cX$-module $\wt\cM\otimes\ccE^{-t/\hb}$ (the definition of the $\cR_\cX$-structure is as above). In particular, $\FcM$ is localized at $\infty$ and $\FcM_\loc$ is the localization of $\FcM$ at $P\moins\{\infty\}$. We know (\cf \cite[Prop.~8.3.1(i)]{Bibi01}) that, under some condition on $\cM$ near $\infty$, the $\cR_\cX$-module $\FcM$ is coherent. Such a condition is satisfied when $\cM$ comes from a (polarized) regular twistor $\cD$-module on~$X$.

Given a (polarized) regular twistor $\cD$-module $(\cM',\cM'',C)$ on $X$, the definition of the sesquilinear pairing $\Fou C$ on $\FcMS'\otimes_{\cO_{\cX|\bS}}\ov{\FcMS''}$ with values in $\Dbh{X}$ needs some care, as it consists in defining a lifting to distributions on $X$ of the localized distribution (or $C^\infty$ function with moderate growth) $\Fou C_\loc$. In \loccit, one first defines $\cFC$ on the total exponential twist $\cFcMS'\otimes\ov{\cFcMS''}$ (where one does not forget the variable~$\tau$); the module $\FcM$ is viewed as the specialization at $\tau=1$ of $\cFcM$; then $\Fou C$ is defined as the specialization (by means of Mellin transform) of $\cFC$.

\subsection{Restriction to $\hb=\hb_o$}
Let us analyze the behaviour of the previous constructions with respect to restriction to $\hb=\hb_o$.

The restriction to $\hb=\hb_o$ of $(\cH^\an,\pi^*\Fou h,\Fou\ccD_\hb)$ is the bundle $H$ equipped with the metric $\Fou h$ and the $\hb_o$-connection $\Fou\ccD_{\hb_o}$. The isomorphism \eqref{eq:isoFF} specializes to an isomorphism
\begin{equation}\label{eq:isoFFrestr}
(H,\Fou h,\Fou\ccD_{\hb_o})\To{{}\cdot e^{(1-\hb_o)\ov t}}(H,e^{2\reel(\hb_o\ov t)}h,\ccD_{\hb_o}-dt)
\end{equation}
and \eqref{eq:isoFFR} specializes to
\begin{equation}\label{eq:isoFFRrestr}
(H,\Fou h,\Fou\ccD_{\hb_o})\To{{}\cdot e^{\ov t-2i\im(\hb_o\ov t)}}(H,h,\ccD_{\hb_o}-(1+\module{\hb_o}^2)dt).
\end{equation}

On the other hand, by flatness of $\cO_\cX[*(P\times\CC)]$ (or $\cO_\cX[*(\{\infty\}\times\CC)]$) over $\cO_\cX$, if $\cM$ is a \emph{strict} $\cR_\cX$-module (that is, it has no $\cO_\CC$-torsion), then so is its localization $\cM_\loc$ or $\wt\cM$. If we put $\ccM_{\hb_o}=\cM/(\hb-\hb_o)\cM$, then the localization ${}_\loc$ or $\wt{\mbox{ }}$ of $\ccM_{\hb_o}$ is the restriction to $\hb=\hb_o$ of the corresponding localization of $\cM$.

Define $\FccM_{\hb_o}$ as $\wt\ccM_{\hb_o}\otimes\cE^{-t/\hb_o}$ if $\hb_o\neq0$ (\ie we twist the $\hb_o$-connection by adding $-dt$) and as $\wt\ccM_0$ with Higgs field obtained by adding $-dt$ if $\hb_o=0$. Then, if $\cM$ is strict, we have $\FccM_{\hb_o}=\FcM/(\hb-\hb_o)\FcM$ and $\FccM_{\hb_o,\loc}=\FcM_\loc/(\hb-\hb_o)\FcM_\loc$.

\section{Proof of Theorem~\ref{th:Fmain}}\label{sec:firstproof}

Let $(\cT,\cS)$ be a polarized regular twistor $\cD$-module of weight $0$ on $\PP^1$ (\ie an object of $\MTr(\PP^1,w)^\rp$, \cf \cite{Bibi01}). We can assume that it takes the form $\cT=(\cM,\cM,C)$ and $\cS=(\id,\id)$. Its restriction to $X^*$ corresponds to a harmonic bundle $(H,h,D_V)$. Taking notation of \eqref{eq:dhat}, we will prove in this section:

\begin{proposition}\label{prop:DRFcM}
The complex $\bR\Gamma(X,\DR\FcM)$ has cohomology in degree $0$ at most, and its nonzero cohomology is a locally free $\cO_{\CC}$-module of finite rank $\wh d$.
\end{proposition}

\subsection{Proof of Theorem \ref{th:Fmain}}

Recall (\cf \cite[Chap.~8]{Bibi01}) that we set $\wt\cM=\cO_{\cX}(*\infty)\otimes_{\cO_{\cX}}\cM$ and, if $p:X\times\wh X\times\CC\to X\times\CC=\cX$ and $\wh p:X\times\wh X\times\CC\to\wh X\times\CC=\wh\cX$ denote the projections, and if ${}\otimes\ccE^{-t\tau/\hb}$ denotes the exponential twist of the $\cR$-structure, we define
\[
\wh\cM\defin\wh p_+p^+(\wt\cM\otimes\ccE^{-t\tau/\hb})=\wh p_+^0p^+(\wt\cM\otimes\ccE^{-t\tau/\hb})\defin\wh p_+^0\cFcM.
\]
The sesquilinear pairing $\cFC$ on $\cFcMS\otimes\ov{\cFcMS}$ is defined in \loccit, and we put $\wh C=\wh p_+^0\cFC$.

\begin{enumerate}
\item
It follows from \cite[Th.~8.4.1]{Bibi01} that, along $\tau=0$, all the necessary conditions for being a polarized regular twistor $\cD$-module (\cf \cite[Def.~4.1.2]{Bibi01}) are satisfied.
\item
Now, the main question concerns the behaviour of $\wh\cM$ out of $\tau=0$. Fix $\tau_o\neq0$ in $\wh X$. Recall (\cite[Prop.~8.3.1(i)]{Bibi01}) that $\wt\cM\otimes\ccE^{-t\tau_o/\hb}$ is $\cR_\cX$-good. It will be enough, in fact, to take $\tau_0=1$, using an obvious homogeneity argument. We denote by $\FcM$ the $\cR_\cX$-module $\wt\cM\otimes\ccE^{-t/\hb}$.

\item
From Proposition \ref{prop:DRFcM} (which holds for any $\tau_o\neq0$), as $\cFcM$ is regular and strictly specializable along $\tau=\tau_o$ and according to \cite{Bibi01} (Proposition~8.3.1(ii) and (iii), Theorem~3.1.8 and \S3.1.d), we obtain:

\begin{corollaire}
For any $\tau_o\neq0$, $\wh\cM$ is strictly specializable and regular along $\tau=\tau_o$ and we have, for any $\alpha\in\CC$,
\[
\psi_{\tau-\tau_o,\alpha}\wh\cM=
\begin{cases}
0&\text{if } \alpha\not\in -\NN^*,\\
\bR^0\Gamma(X,\DR\wt\cM\otimes\ccE^{-t\tau_o/\hb})&\text{if }\alpha\in-\NN^*.
\end{cases}
\]
\end{corollaire}

\item
This corollary implies that, near any $\tau_o\neq0$, $\wh\cM$ is equal to the level $-1$ of its $V$-filtration along $\tau=\tau_o$. By regularity, $\wh\cM$ is therefore $\cO_{\wh\cX}$-coherent and, as $\dim\wh\cM/(\tau-\tau_o)\wh\cM=\dim\psi_{\tau-\tau_o,\alpha}\wh\cM=\wh d$ is independent of $\tau_o\neq0$, $\wh\cM$ is $\cO_{\wh\cX}$-locally free of rank $\wh d$ out of $\tau=0$. Its characteristic variety in $T^*(\wh X\moins\{0\})\times\CC$ is equal to $\text{zero section}\times\CC$, and its characteristic variety in $T^*(\wh X)\times\CC$ is contained in
\[
\big(\text{zero section}\cup T^*_0\wh X\big)\times\CC,
\]
so $\wh\cM$ is holonomic (\cf \cite[Def.~1.2.4]{Bibi01}). The corollary also implies that S-decom\-po\-sability is trivially satisfied near $\tau_o\neq0$. We have obtained therefore Condition (HSD) of \cite[Def.~4.1.2]{Bibi01}.

\item
At this step, we know that $\wh \cM$ is $\cO_{\wh\cX}$-locally free of finite rank out of $\tau=0$. It follows (\cf \cite[Lemma~1.5.3]{Bibi01}) that, on this domain, $(\wh\cM,\wh\cM,\wh C)$ is a smooth object of $\RTriples(\wh X^*)$. We wish to show that $\wh C$ defines, by gluing, a family parametrized by $\wh X^*$ of trivial vector bundles on $\PP^1$, coming from a $C^\infty$ vector bundle $\wh H$ on $\wh X^*$ equipped with a Hermitian metric $\wh h$ by  the correspondence of \cite[Lemma~2.2.2]{Bibi01}. By this construction, the metric will be \emph{harmonic}. By a simple homogeneity argument with respect to $\tau$, it is enough to prove this property in some neighbourhood of $\tau=0$, that we still denote by $\wh X^*$.

As we know, by \cite[Theorem~8.4.1]{Bibi01}, that the twistor properties are satisfied along $\tau=0$ by $\wh\cT=(\wh\cM,\wh\cM,\wh C)$ equipped with the polarization $\wh\cS=(\id,\id)$, we can apply the argument given in \cite[\oldS\S5.4.c--5.4.e]{Bibi01} to get the twistor property and the polarizability in some neighbourhood of $\tau=0$.\qed
\end{enumerate}

\medskip
We will now prove Proposition \ref{prop:DRFcM}. As $\FcM$ is a good $\cR_\cX$-module, we \emph{a priori} know that the cohomology of $\bR\Gamma(X,\DR\FcM)$ is $\cO_\CC$-coherent. It is then enough to prove that, for any $\hb_o\in\CC$, the complex $\bR\Gamma(X,\DR\FccM_{\hb_o})$ has cohomology in degree~$0$ only, and that the dimension of $\bH^0(X,\DR\FccM_{\hb_o})$ is $\wh d$ (let us recall that $\FccM_{\hb_o}=\FcM/(\hb-\hb_o)\FcM$).

As in \cite{Zucker79}, we will identify the complex $\DR\FccM_{\hb_o}$ with a $L^2$ complex. This identification is \emph{local} on $X$. The $L^2$ cohomology on $X$ is then obtained by $L^2$-Hodge theory. That the dimension of $\bH^\bbullet(X,\DR\FccM_{\hb_o})$ is independent of $\hb_o$ will follow from the independence of the corresponding Laplacian with respect to $\hb_o$ (an argument which comes from \cite{Simpson92}).

\subsection{The meromorphic $L^2$ de~Rham and Dolbeault complexes}
In order to give a common proof whether $\hb_o=0$ or not, it will be convenient to consider the twisted module $\FccM_{\hb_o}\otimes\ccE^{-c(\hb_o)t}$, where $c(\hb_o)$ denotes the usual conjugate of $\hb_o$, so that $\module{\hb_o}^2=\hb_oc(\hb_o)$ (we keep the more traditional notation $\ov\hb_o$ for the ``geometric conjugate'' $-1/\hb_o$). In other words, $\FccM_{\hb_o}\otimes\ccE^{-c(\hb_o)t}$ is nothing but $\wt\ccM_{\hb_o}$ as a $\cO_\cX$-module, equipped with the twisted $\hb_o$-connection $\ccD'_{\hb_o}-(1+\module{\hb_o}^2)dt$.

Recall that we denote by $\FccM_{\loc,\hb_o}$ the localized module of $\FccM_{\hb_o}$ at all points of~$P$ (but $\infty$ is unnecessary as $\FccM_{\hb_o}$ is yet localized at $\infty$). We consider the meromorphic $L^2$ complex $\DR(\FccM_{\loc,\hb_o}\otimes\ccE^{-c(\hb_o)t})_{(2)}$ obtained by taking sections of $\ccM_{\loc,\hb_o}$ or $\ccM_{\loc,\hb_o}\otimes\Omega_X^1$ which are locally $L^2$, as well as their image by $\ccD'_{\hb_o}-(1+\module{\hb_o}^2)dt$, when one takes the metric $h$ on the restriction $V_{\hb_o}$ of $\ccM_{\loc,\hb_o}$ to $X^{*\an}$ ($V_{\hb_o}$ is the holomorphic subbundle of $H$ defined by the $d''$ operator $\ccD''_{\hb_o}=D''_V+(\hb_o-1)\theta''_E$) and a metric locally equivalent to the Poincar\'e metric near each puncture $P$ on $X^*$. We have a natural morphism
\[
\DR(\FccM_{\loc,\hb_o}\otimes\ccE^{-c(\hb_o)t})_{(2)}\to \DR(\FccM_{\hb_o}\otimes\ccE^{-c(\hb_o)t}).
\]
Indeed, out of $\infty$, this was explained in \cite[\S6.2.a]{Bibi01} (as the terms of the left-hand complex are not evidently contained in the corresponding terms of the right-hand one). Near $\infty$, the inclusion is clear, as $\FccM_{\hb_o}$ is equal there to $\FccM_{\loc,\hb_o}$.

\begin{lemme}\label{lem:FL2mero}
The natural morphism $\DR(\FccM_{\loc,\hb_o}\otimes\ccE^{-c(\hb_o)t})_{(2)}\to \DR(\FccM_{\hb_o}\otimes\ccE^{-c(\hb_o)t})$ is a quasi-isomorphism.
\end{lemme}

\begin{proof}
This was proved in \loccit, Prop.~6.2.4, out of $\infty$. We therefore consider the situation near $\infty$, with local coordinate $t'$ and we forget the index ``$\loc$'', as $\FccM_{\hb_o}$ is equal to its localized module near $t'=0$.

By the regularity assumption of $\cM$ near $\infty$, we know that there exists a local meromorphic basis $\bme^{(\hb_o)}$ of $\wt\ccM_{\hb_o}$ in which the matrix of $\ccD'_{\hb_o}$ has a simple pole at $t'=0$ (\cf \cite[Formula (5.3.7)]{Bibi01}). By considering the maximal order of the poles of the coefficients of a section of $\wt\ccM_{\hb_o}$ on the basis $\bme^{(\hb_o)}$, and according to the term \hbox{$(1+\module{\hb_o}^2)dt'/t^{\prime2}$} in the $\hb_o$-connection, one obtains that \hbox{$\cH^{-1}\big(\DR(\FccM_{\hb_o}\otimes\ccE^{-c(\hb_o)t})\big)=0$}, hence also $\cH^{-1}\big(\DR(\FccM_{\hb_o}\otimes\ccE^{-c(\hb_o)t})_{(2)}\big)=0$.

On the other hand, by the same reasoning, any local section at $t'=0$ of $\wt\ccM_{\hb_o}\otimes\Omega^1_X$ with maximum order of pole equal to $k$ is equivalent, modulo the image of $\ccD'_{\hb_o}+(1+\module{\hb_o}^2)dt'/t^{\prime2}$, to a section having a pole of maximum order $\leq k-1$. Iterating this process, and according to the moderate behaviour of the $h$-norm of each element of the basis $\bme^{(\hb_o)}$, we see that such a section is equivalent to a section of $\wt\ccM_{\hb_o}\otimes\Omega^1_X$ which is $L ^2$ with respect to the metric $h$, or equivalently that
$$
\cH^0\big(\DR(\FccM_{\hb_o}\otimes\ccE^{-c(\hb_o)t})_{(2)}\big)\to
\cH^0\big(\DR(\FccM_{\hb_o}\otimes\ccE^{-c(\hb_o)t})\big)
$$
is onto.

Last, given a local section of $\wt\ccM_{\hb_o}\otimes\Omega^1_\cX$ which which $L^2$ (with respect to the metric~$h$) and is in $\big(\ccD'_{\hb_o}-(1+\module{\hb_o}^2)\big)\wt\ccM_{\hb_o}$, the same kind of argument shows that it is in the image of  a $L^2$ section of $\wt\ccM_{\hb_o}$; equivalently,
$$
\cH^0\big(\DR(\FccM_{\loc,\hb_o}\otimes\ccE^{-c(\hb_o)t})_{(2)}\big)\to
\cH^0\big(\DR(\FccM_{\loc,\hb_o}\otimes\ccE^{-c(\hb_o)t})\big)
$$
is injective.
\end{proof}

\subsection{The $L^2$ de~Rham-Dolbeault Lemma}
We consider the $C^\infty$ bundle $H$, equipped with the metric $h$ and with the $\hb_o$-connection $\ccD_{\hb_o}-(1+\module{\hb_o}^2)dt$, that we denote below by $\wt\ccD_{\hb_o}$ for simplicity, and the associated $L^2$ complex $\cLd^\cbbullet(H,h,\wt\ccD_{\hb_o})$. In particular, one can notice that the $d''$-operator is $\ccD''_{\hb_o}$, the corresponding holomorphic subbundle is $V_{\hb_o}$, and the extension of this holomorphic subbundle obtained by considering sections, the $h$-norm of which has moderate growth, is $\ccM_{\loc,\hb_o}$ (\cf \cite[Cor.~5.3.1(1)]{Bibi01}).

The ``holomorphic'' $L^2$ subcomplex is the following subcomplex of the $L^2$ complex $\cLd^\cbbullet(H,h,\wt\ccD_{\hb_o})$:
\begin{equation}\label{eq:FholL2}
0\to\ker\ccD_{\hb_o}^{\prime\prime(0)}\To{\wt\ccD'_{\hb_o}}\ker\ccD_{\hb_o}^{\prime\prime(1)}\cap\cLd^{(1,0)}(H,h,\ccD''_{\hb_o}) \to0,
\end{equation}
where we denote by $\ccD_{\hb_o}^{\prime\prime(k)}$ the action of $\ccD''_{\hb_o}$ on $\cLd^k(H,h,\ccD''_{\hb_o})$. Our aim in this paragraph is to prove:

\begin{lemme}[$L^2$ de~Rham-Dolbeault Lemma]\label{lem:FDolb}
Assume that $(H,h,D_V)$ is a tame harmonic bundle on $X^{*\an}$. Then, the inclusion of the holomorphic $L^2$ subcomplex \eqref{eq:FholL2} in $\cLd^\cbbullet(H,h,\wt\ccD_{\hb_o})$ is a quasi-isomorphism.
\end{lemme}

The proof will be analogous to the proof of the Dolbeault Lemma in \cite{Zucker79} and will be parallel to that of Theorem 6.2.5 in \cite{Bibi01}, which we constantly refer to. As above, we will work near $\infty$, as the result out of $\infty$ is given by \oldS\S6.2.d and 6.2.e of \loccit

In the definition of the $L^2$ complex, the $L^2$ condition on sections is the same as in \loccit, as well as the action of the anti-holomorphic part of the connection. The $L^2$ condition on the derivative of sections is changed. The new term $(1+\module{\hb_o}^2)dt'/t^{\prime2}$ in the holomorphic part of the connection will make proofs easier. 

We use polar coordinates: $t'=re^{i\theta}$. Let us first recall some notation of \loccit Near $t'=0$, the bundle $H$ is equipped with a $\ccD''_{\hb_o}$-holomorphic basis $\bme^{\prime(\hb_o)}$. The $h$-norm of the elements in this basis has moderate growth near $t'=0$. We denote them by $e^{\prime(\hb_o)}_{\beta,\ell,k}$, where $\beta=\beta'+i\beta''$ runs in a finite set of complex numbers with real part $\beta'$ in $[0,1[$, $\ell$ is an integer (the weight of the element) and $k$ is an index used to distinguish the various elements having the same $\beta$ and $\ell$. In such a basis, the matrix of $\ccD'_{\hb_o}$ is written $\Theta'_{\hb_o}$. It is decomposed as $\Theta'_{\hb_o,\diag}+\Theta'_{\hb_o,\nilp}$, with
\begin{align*}
\Theta'_{\hb_o,\diag}&=\oplus_\beta(q_{\beta,\imhb_o}+\beta)\star\hb\id\,\dfrac{dt'}{t'}\\
\Theta'_{\hb_o,\nilp}&=\big[\rY+P(t,\hb)\big]\dfrac{dt'}{t'},
\end{align*}
with $\rY=(\oplus_\beta \rY_\beta)$, and $q_{\beta,\imhb_o}$ is an integer chosen such that, denoting by $\imhb_o$ the imaginary part of $\hb_o$, the number $\ell_{\hb_o}(q_{\beta,\imhb_o}+\beta)\defin q_{\beta,\imhb_o}+\beta'-\imhb_o\beta''$ belongs to $[0,1[$. The numbering of the basis is such that we have, for any $\ell,k$, $\rY(e^{\prime(\hb_o)}_{\beta,\ell,k})=e^{\prime(\hb_o)}_{\beta,\ell-2,k}$ and a formula for $P(t,\hb_o)$ is given by (6.2.7) in \loccit

Recall that we denote $\wt\ccD_{\hb_o}=\ccD_{\hb_o}+(1+\module{\hb_o}^2)dt'/t^{\prime2}$. Then, with obvious notation, $\wt\Theta'_{\hb_o,\nilp}=\Theta'_{\hb_o,\nilp}$ and $\wt\Theta'_{\hb_o,\diag}=\Theta'_{\hb_o,\diag}+(1+\module{\hb_o}^2)\id dt'/t^{\prime2}$.

\begin{proof}[Vanishing of $H^2$]
First, Lemma 6.2.11 of \cite{Bibi01}, with a fixed value $\hb=\hb_o$, applies here without any modification. We are therefore led to proving that, if for any $\beta$ with $\ell_{\hb_o}(q_{\beta,\imhb_o}+\beta)=0$ and any $\ell\leq-1$ (in fact, $\ell=-1$ would suffice as $\hb$ is fixed to $\hb_o$ here), $f(r)e'_{\beta,\ell,k}\, \frac{dt'}{t'}\wedge\frac{d\ov{t'}}{\ov{t'}}$ is a local section of $\cLd^2(H)$, then it belongs to the image of~$\wt\ccD_{\hb_o}$.

Remark that
\begin{multline*}
\wt\ccD_{\hb_o}\Big(t'f(r)e'_{\beta,\ell,k}\Big(\hb_o\frac{dt'}{t'}+\frac{d\ov{t'}}{\ov{t'}}\Big)\Big)\\
=\big(1+\module{\hb_o}^2+\hb_o+(\beta\star\hb_o)t'\big)f(r)e'_{\beta,\ell,k}\, \frac{dt'}{t'}\wedge\frac{d\ov{t'}}{\ov{t'}}\\
+\Theta'_{\hb_o,\nilp}\Big(t'f(r)e'_{\beta,\ell,k}\Big(\hb_o\frac{dt'}{t'}+\frac{d\ov{t'}}{\ov{t'}}\Big)\Big).
\end{multline*}
As in \loccit, the last term is easily seen to be $L^2$ and, because of the factor $t'$, is in the image of $\wt\ccD''_{\hb_o}=\ccD''_{\hb_o}$ by Lemma 6.2.11 of \cite{Bibi01}. For the same reason, the part which is multiple of $t'$ in the middle term is in the image of $\ccD''_{\hb_o}$. Both are therefore in the image of $\wt\ccD_{\hb_o}$. To conclude, it remains to notice that the constant $1+\module{\hb_o}^2+\hb_o$ cannot be equal to $0$.
\end{proof}

\begin{proof}[Computation of $H^1$]
From the previous result, the $L^2$-complex $\cLd^\cbbullet(H,h,\wt\ccD_{\hb_o})$ is quasi-isomorphic to its subcomplex
\[
0\to\cLd^0(H,h,\wt\ccD_{\hb_o})\To{\wt\ccD_{\hb_o}}\ker\wt\ccD^{(1)}_{\hb_o}\to0.
\]
We will now show the analogue of Lemma 6.2.13 in \cite{Bibi01}. That is, we will show that any local section $\psi\,dt/t+\varphi\,d\ov t/\ov t$ of $\ker\wt\ccD^{(1)}_{\hb_o}\subset\cLd^1(H,h,\wt\ccD_{\hb_o})$ can be written as the sum of a term in $\im\wt\ccD_{\hb_o}$ and a term in $\cLd^{(1,0)}(H,h)\cap\ker\wt\ccD_{\hb_o}^{(1)}$.

The first\label{pageerr} part of the proof of Lemma 6.2.13 in \loccit applies similarly to the present situation, and reduces the proof to the case where we start with a local section $\omega=\psi\sfrac{dt'}{t'}+\varphi\sfrac{d\ov{t'}}{\ov{t'}}$ of $\ker\wt\ccD^{(1)}_{\hb_o}$, with $\varphi=\sum_{\beta,\ell,k}\varphi_{\beta,\ell,k}(r)e^{\prime(\hb_o)}_{\beta,\ell,k}$ and satisfying $\wt\ccD_{\hb_o}\omega=0$.

We consider then the component on $e^{-i\theta} e^{\prime(\hb_o)}_{\beta,\ell,k}\,\frac{dt'}{t'}\wedge\frac{d\ov{t'}}{\ov{t'}}$ of the relation $\ccD''_{\hb_o}(\psi dt'/t')+\wt\ccD'_{\hb_o}(\varphi d\ov{t'}/\ov{t'})=0$. Denoting by $\psi_{\beta,\ell,k;-1}(r)$ the component on $e^{-i\theta}$ in the Fourier expansion of $\psi_{\beta,\ell,k}$, we get, for any $\beta,\ell,k$,
\begin{align*}
\varphi_{\beta,\ell,k}(r)e^{\prime(\hb_o)}_{\beta,\ell,k}\frac{d\ov{t'}}{\ov{t'}}&=\frac12 r(r\partial_r-1)\psi_{\beta,\ell,k;-1}(r)e^{\prime(\hb_o)}_{\beta,\ell,k}\frac{d\ov{t'}}{\ov{t'}}\\
&=\ccD''_{\hb_o}\big(re^{-2i\theta}\psi_{\beta,\ell,k;-1}(r) e^{\prime(\hb_o)}_{\beta,\ell,k}\big).
\end{align*}
As $\psi dt'/t'$ is $L^2$, it follows that $r\psi$ is $L^2$, and that $\wt\ccD'_0\big(re^{-2i\theta}\psi_{\beta,\ell,k;-1}e^{\prime(\hb_o)}_{\beta,\ell,k})$ is also~$L^2$.

This computation shows that $\omega$ is equivalent, modulo $\im\wt\ccD_{\hb_o}$, to a $(1,0)$-section which is $L^2$, and which is in $\ker\wt\ccD^{(1)}_{\hb_o}$ (because $\wt\ccD_{\hb_o}\omega=0$), as wanted.
\end{proof}

\pagebreak[2]
\subsection{End of the proof of Proposition \ref{prop:DRFcM}}
We argue in four steps.
\begin{enumerate}
\item
Arguing exactly as in \cite[\S6.2.f]{Bibi01}, we show that the ``holomorphic'' $L^2$ complex \eqref{eq:FholL2} is equal to its subcomplex $\DR(\FccM_{\loc,\hb_o}\otimes\ccE^{-c(\hb_o)t})_{(2)}$. By coherence, the hypercohomology of $\DR(\FccM_{\hb_o}\otimes\ccE^{-c(\hb_o)t})$ is finite dimensional. By Lemma \ref{lem:FL2mero} and by the argument above, so is the hypercohomology of the holomorphic $L^2$ complex \eqref{eq:FholL2}.
\item
From Lemma \ref{lem:FDolb} and the previous result, we conclude that the cohomology of $\Gamma\big(X,\cLd^\cbbullet(H,h,\wt\ccD_{\hb_o})\big)$ is finite dimensional. According to the isometry \eqref{eq:isoFFRrestr}, the cohomology of $\Gamma\big(X,\cLd^\cbbullet(H,\Fou h,\Fou\ccD_{\hb_o})\big)$ is finite dimensional. We can therefore apply Hodge Theory to this $L^2$ cohomology. The corresponding space of harmonic $k$-forms ($k=0,1,2$) is finite dimensional, and its dimension does not depend on $\hb_o$, as the Laplacian of $\Fou\ccD_{\hb_o}$ with respect to the metric $\Fou h$ is essentially independent of $\hb_o$, due to the harmonicity property of the triple $(H,\Fou h,\Fou D_V)$.
\item
Going in the reverse direction, we see that the dimension of the space $\bH^k\big(X,\DR(\FccM_{\hb_o}\otimes\ccE^{-c(\hb_o)t})\big)$ ($k=-1,0,1$) is independent of $\hb_o$. If $\hb_o=1$, there is cohomology in degree $0$ only (this is well-known for a regular holonomic $\cD_X$-module twisted by an exponential $e^{\lambda t}$, with $\lambda\in\CC^*$). This is therefore true for any $\hb_o$, and moreover the dimension of $\bH^0$ is independent of $\hb_o$.
\item
It remains to remark that the hypercohomologies of $\DR(\FccM_{\hb_o}\otimes\ccE^{-c(\hb_o)t})$ and $\DR(\FccM_{\hb_o})$ have the same dimension. This is clear if $\hb_o=0$ as both objects are equal. On the other hand, if $\hb_o\neq0$, we are reduced to a question on $\cD_X$-modules. Working algebraically, we are reduced to proving that, given a regular holonomic $\CC[t]\langle\partial_t\rangle$-module, the dimension of the cokernel of
\[
\partial_t-\lambda:M\to M
\]
does not depend on $\lambda\in\CC^*$. This is a follows from the regularity of $M$ at infinity.\qed
\end{enumerate}

\begin{remarque}
It would have been possible to give a different proof of Proposition \ref{prop:DRFcM} when $\hb_o\neq0$, by using the $\hb_o$-connection $\ccD_{\hb_o}-dt$ on $H$ with the metric $e^{2\reel(\hb_o\ov t)}h$. This proof would be analogous to that of \cite{Bibi98}\footnote{Let us take this opportunity to correct a minor mistake in \cite{Bibi98}: on page 1283, the inequality is
\[
\int_\rho^{r_1}r^{2\beta}\module{\log r}^{k-2}\psi(r)\,\frac{dr}{r}\leq \rho^{2\beta}\module{\log \rho}^{k-2}\psi(\rho)(\module{\log\rho}-\module{\log r_1})
\]
and the constant $C$ is bounded by $4\module{\log r_1}^{-1}<+\infty$. Similarly, on page 1284, $\ell$.~5, the constant $C$ is bounded by $4\kappa(\epsilon)\module{\log r_1}^{-1}<+\infty$.} and we would use the isometry \eqref{eq:isoFFrestr} instead of \eqref{eq:isoFFRrestr}. Nevertheless, the intermediate steps would be different, as the analogue of Lemma \ref{lem:FL2mero}, where the $L^2$ condition is taken with respect to the metric $e^{2\reel(\hb_o\ov t)}h$, is not true. One has to work, as in \cite{Bibi98}, in the space obtained from $X$ by a real blowing-up at infinity. The comparison between various complexes has to be made on this space. However, such a proof does not seem to extend to $\hb_o=0$. This is why it is not developed here.
\end{remarque}

\providecommand{\bysame}{\leavevmode ---\ }
\providecommand{\og}{``}
\providecommand{\fg}{''}
\providecommand{\smfandname}{\&}
\providecommand{\smfedsname}{eds.}
\providecommand{\smfedname}{ed.}

\part*{Erratum to ``Fourier-Laplace transform of irreducible regular differential systems on the Riemann sphere''}
\setcounter{section}{0}
\setcounter{subsection}{0}
\setcounter{lemme}{0}
\def\thefootnote{*}

There\footnote{This erratum was written in may 2007.}
 are two mistakes in the proof of Theorem~1 above:
\begin{enumerate}
\item\label{err:1}
In \S3.1, point 5, we assert:

``By simple homogeneity considerations with respect to $\tau$, it suffices to prove the property in the neighbourhood of $\tau=0$.''

It happens that homogeneity does not lead to such a statement. One has to prove the twistor property for the pairing $\wh C$ at any $\tau^o\neq0$. While the proof given above holds for $|\tau^o|$ small enough by an argument of degeneration, we will give in \S\ref{secapp:B} a proof for any fixed $\tau^o$ and it is enough for such a proof to give the argument when $\tau^o=1$.

\item\label{err:2}
In the proof of Lemma 4 (a main tool for Proposition~1), the computation of~$H^1$ cannot follow the same lines as Lemma 6.2.13 in \cite{errBibi01c}, since this lemma contained a mistake (which has can easily be corrected in the tame case of \cite{errBibi01c}, but not in the present context). We will instead use the argument indicated in Remark~1 above.
\end{enumerate}

In this erratum, we correct these two points. The correction for \eqref{err:2} is given in \S\ref{secapp:A} and that for \eqref{err:1} in \S\ref{secapp:B}.

\subsubsection*{Acknowledgements}
I~gratefully thank the referee of \cite{errBibi05} for having pointed out these mistakes and for having given a suggestion for their correction. In particular, Lemma~\ref{lem:decay} is due to him.

\section{Correction of the proof of Proposition~1}\label{secapp:A}
The proof of Proposition~1 follows the same lines as in \S3.4 above once we have proved the lemma below. Nevertheless, instead of using the isometry (2.5), we will use (2.4) of \loccit\ In order to simplify the notation, we will set in the following $\wt\ccD_\hb=\ccD_\hb-dt$ (take care that this does not correspond to the notation above) and $h_\hb=e^{2\reel\hb\ov t}\pi^*h$, in particular, $h_\hbo=e^{2\reel\hbo\ov t}h$.

\begin{lemme}\label{lem:comparison}
For any $\hbo\in\Omega_0$, there is an isomorphism in the derived category $D^b(\CC_{\PP^1})$:
\[
\DR(\FccM_\hbo)\simeq \cLd^{1+\cbbullet}(H,h_\hbo,\wt\ccD_\hbo).
\]
\end{lemme}

Let us quickly recall the notation used above, after \cite{errBibi01c}. We denote by $\Omega_0$ the complex line with coordinate $\hb$ (or an open neighbourhood of the closed disc \hbox{$|\hb|\leq1$}). If $X$ is a complex manifold (here, $X=\PP^1$ or $X$ is a disc), then $\cX=X\times\Omega_0$ (\eg $\cP^1=\PP^1\times\Omega_0$) and $\pi:\cX\to X$ denotes the projection. We will denote by $p:\cX\to\Omega_0$ the other projection. We also set $\rL(t)=\big|\log|t|^2\big|$.

\begin{proof}[Proof of Lemma \ref{lem:comparison}]
We will distinguish whether $\hbo=0$ or not. When \hbox{$\hbo=0$} we continue using Lemma~4 (Dolbeault Lemma) above as it stands (with the supplementary assumption that $\hbo=0$), corrected as in \S\ref{subsec:Dolbeault} below. It says that the natural inclusion $\DR(\FccM_0)_{(2)}\hto\cLd^{1+\cbbullet}(H,h,\wt\ccD_0)$ is a quasi-isomorphism. We then use Lemma~3 of above to conclude.

When $\hbo\neq0$, we will change the argument, and use that indicated in Remark~1 above. This will be done in \S\ref{subsec:Poincare} below.
\end{proof}

\subsection{The Dolbeault lemma}\label{subsec:Dolbeault}
We correct the statement given on page \pageref{pageerr}, line~4 above (page 1178, line~1 of the published version) when $\hbo=0$. Let $\omega=\psi\dfrac{dt'}{t'}+\varphi\dfrac{d\ov{t'}}{\ov{t'}}$ be in $\ker \wt\ccD_0^{(1)}$. We wish to prove that, modulo the image by $\wt\ccD_0$ of $\cLd^0(H,h,\wt\ccD_0)$, we can reduce $\varphi$ to be as written in page \pageref{pageerr}, line~7 above (page 1178, line~3 of the published version). Expanding $\omega$ on the basis $(e_{\beta,\ell,k}^{\prime(0)})$, the $L^2$ condition reads: $\mtp^{\beta'}\Ltp^{\ell/2}\module{\psi_{\beta,\ell,k}}$ and $\mtp^{\beta'}\Ltp^{\ell/2}\module{\varphi_{\beta,\ell,k}}$ belong to $L^2(d\theta\,dr/r)$. Then $\wt\ccD_0^{(1)}\omega=0$ reads, setting $\beta=\beta'+i\beta''$ and $\alpha(t')=1/(1+i\beta''t'/2)$,
\[
-\ov{t'}\partial_{\ov{t'}}\psi_{\beta,\ell,k}+\frac{1}{\alpha(t')t'}\varphi_{\beta,\ell,k}+\xi_{\beta,\ell,k}=0,
\]
with $\xi=\sum\xi_{\beta,\ell,k}e_{\beta,\ell,k}^{\prime(0)}$ defined by $\xi \dfrac{dt'}{t'}\wedge \dfrac{d\ov{t'}}{\ov{t'}}=\Theta'_{0,\nilp}\varphi$. It follows that
\[
\varphi_{\beta,\ell,k}=\ov{t'}\partial_{\ov{t'}} (\alpha(t')t'\psi_{\beta,\ell,k})-\alpha(t')t'\xi_{\beta,\ell,k}.
\]
Firstly, $\mtp^{\beta'}\Ltp^{-1+\ell/2}\module{t'\psi_{\beta,\ell,k}}$ clearly belongs to $L^2(d\theta\,dr/r)$, hence on the one hand, $t'\psi_{\beta,\ell,k}e_{\beta,\ell,k}^{\prime(0)}$ is a section of $\cLd^0(H,h)$. On the other hand,
\[
\Theta'_0(\alpha(t')t'\psi_{\beta,\ell,k}e_{\beta,\ell,k}^{\prime(0)})=\psi_{\beta,\ell,k}e_{\beta,\ell,k}^{\prime(0)}\frac{dt'}{t'}+\Theta'_{0,\nilp}(\alpha(t')t'\psi_{\beta,\ell,k}e_{\beta,\ell,k}^{\prime(0)}).
\]
As we know that $\Theta'_{0,\nilp}$ is bounded with respect to the $L^2$ norms, it follows that the left-hand term is a section of $\cLd^{(1,0)}(H,h)$. We conclude that $\wt\ccD_0(\alpha(t')t'\psi)$ is~$L^2$ and that the $(0,1)$-part of $\omega-\wt\ccD_0(\alpha(t')t'\psi)$ is equal to $-\alpha(t')t'\xi \sfrac{d\ov{t'}}{\ov{t'}}$.

Secondly, by the property of $\Theta'_{0,\nilp}$, we find that $\mtp^{\beta'}\Ltp^{1+\ell/2}\module{\xi_{\beta,\ell,k}}$ also belongs to $L^2(d\theta\,dr/r)$. Let us now argue as in \cite[Lemma 6.2.11]{errBibi01c}. We expand $\xi_{\beta,\ell,k}$ as a Fourier series $\sum_n\xi_{\beta,\ell,k,n}(r)e^{in\theta}$ with $r=\mtp$, and set $\xi_{\beta,\ell,k,\neq0}=\xi_{\beta,\ell,k}-\xi_{\beta,\ell,k,0}$. We then find that it is possible to solve $\ov{t'}\partial_{\ov{t'}}\eta_{\beta,\ell,k,\neq0}=\xi_{\beta,\ell,k,\neq0}$ with $\eta_{\beta,\ell,k,\neq0}$ being a local section of $\cLd^0(H,h)$. As above, we then show that $\Theta'_0(\alpha(t')t'\eta_{\beta,\ell,k,\neq0}e_{\beta,\ell,k}^{\prime(0)})$ is a section of $\cLd^{(1,0)}(H,h)$.

We finally conclude that $\omega-\wt\ccD_0[\alpha(t')t'(\psi-\eta_{\neq0})]$ satisfies the desired property.

\subsection{The Poincar\'e lemma}\label{subsec:Poincare}
We will now give the proof of Lemma~\ref{lem:comparison} when $\hbo\neq0$, a condition that we assume to hold for the remaining of this subsection.

\subsubsection*{Reduction of the proof of Lemma~\ref{lem:comparison} to local statements when $\hbo\neq0$}
We will first work with the metric $h$ (and not $h_\hbo$). We denote by $\FccM_{\hbo,\loc}$ the localization of $\FccM_\hbo$ at the singularities $P$ (note that, at infinity, $\FccM_\hbo$ is yet equal to its localized module) and by $\DR(\FccM_{\hbo,\loc})_{(2),h}$ the meromorphic $L^2$ de~Rham complex, which is a subcomplex of $\DR\FccM_{\hbo,\loc}$. In fact, it is a subcomplex of $\DR\FccM_\hbo$: at finite distance, this is \cite[Prop.~6.2.4]{errBibi01c} and at infinity this is clear. The argument of Lemma~3 above gives:

\begin{lemme}\label{lem:new3}
The inclusion of complexes $\DR(\FccM_{\hbo,\loc})_{(2),h}\hto\DR\FccM_\hbo$ is a quasi-isomorphism.\qed
\end{lemme}

On the other hand, by definition, $\DR(\FccM_{\hbo,\loc})_{(2),h}$ is a sub-complex of the $L^2$ complex $\cLd^{1+\cbbullet}(H,h,\wt\ccD_\hbo)$ and, according to \cite[Th.~6.2.5]{errBibi01c}, the inclusion is a quasi-isomorphism at finite distance. Lemma~\ref{lem:comparison} now follows from the following two statements:

\begin{enumerate}\let\theenumi\theequation
\refstepcounter{equation}
\item\label{enum:DR2h}
The natural inclusion $\DR(\FccM_{\hbo,\loc})_{(2),h}\hto\cLd^{1+\cbbullet}(H,h,\wt\ccD_\hbo)$ is a quasi-isomorphism near $\infty$.
\refstepcounter{equation}
\item\label{enum:hhhbo}
Both inclusions of complexes
\[
\cLd^{1+\cbbullet}(H,h,\wt\ccD_\hbo)\hfrom\cLd^{1+\cbbullet}(H,h+h_\hbo,\wt\ccD_\hbo)\hto\cLd^{1+\cbbullet}(H,h_\hbo,\wt\ccD_\hbo).
\]
are quasi-isomorphisms.
\end{enumerate}

Both questions are now local near $\infty$, and we will restrict to an open disc at $\infty$. So, we set $t'=1/t$ and we denote by $X$ the open disc centered at $0$ and of radius $r_0<1$ in~$\CC$, with coordinate $t'$, and we set $X^*=X\moins\{0\}$. We still keep the notation $h_\hbo$ for the metric $e^{2\reel\hbo/{\ov t'}}h$. We will work with polar coordinates with respect to $t'$.

\subsubsection*{The setting}
We consider the real blow-up
\[
\rho:\wt X\defin [0,r_0[{}\times S^1\to X,\quad (r,\theta)\mto t'=re^{i\theta}.
\]

We will use the sheaf $\cA^{\rmod}_{\wt X}$ on $\wt X$, consisting of holomorphic functions on $\wt{X^*}\!=\!X^*$ which have moderate growth along $r=0$. It is known that $\cA^{\rmod}_{\wt X}$ is stable by~$\partiall_{t'}$. We also consider the differential $1$-forms on $\wt X$:
\begin{align*}
\omega_r&=\frac{(\hbo+1)}{2}\,\frac{dr}{r}+i\frac{(\hbo-1)}{2}\,d\theta\\
\omega_\theta&=-i\frac{(\hbo-1)}{2}\,\frac{dr}{r}+\frac{(\hbo+1)}{2}\,d\theta,
\end{align*}
which form a basis of $1$-forms and which satisfy
\[
\frac{dr}r-id\theta=\omega_r-i\omega_\theta,\quad \frac{dr}r+id\theta=\frac1\hbo\,(\omega_r+i\omega_\theta).
\]
Let us denote by $d$ the differential. The decomposition $d=d'+d''$ on $X$ can be lifted to $\wt X$ and, for a $C^\infty$ function $\varphi(r,\theta)$ on $\wt X$, we have
\[
(d''+\hbo d')\varphi= r\partial_r(\varphi)\omega_r+\partial_\theta(\varphi)\omega_\theta.
\]
Similarly, for a $1$-form $\eta=\varphi\omega_r+\psi\omega_\theta$, we have
\[
(d''+\hbo d')\eta= \bigl(r\partial_r(\psi)-\partial_\theta(\varphi)\bigr)\omega_r\wedge\omega_\theta.
\]

\subsubsection*{The $L^2$ complexes}
Recall that, in this local setting, we denote by $\wt\ccD_\hbo$ the connection $\ccD_\hbo+dt'/t^{\prime2}$. We are interested in computing the cohomology of the complex $\cLd^{1+\cbbullet}(H,\varh,\wt\ccD_\hbo)$, where $\varh$ denotes one of the metrics $h$, $h_\hbo$ or $h+h_\hbo$, which is defined exactly like in \cite[\S6.2.b]{errBibi01c}.

We can similarly define the corresponding $L^2$ complex $\wcLd^{1+\cbbullet}$ by working on~$\wt X$. Let us notice that the use of polar coordinates is convenient to express the $L^2$ condition.

The local basis $\bme^{\prime(\hbo)}\defin(e^{\prime(\hbo)}_{\beta,\ell,k})$ which was introduced in \cite{errBibi01c} for the bundle $(H,h,\ccD_\hbo)$ remains holomorphic with respect to $\wt\ccD''_\hbo$, and also $L^2$-adapted for the metric $\varh$ (in \loccit, we used the notation $\bme^{\prime(\hbo)}$ for a frame defined when $\hb$ varies; in this paragraph, we reduce it modulo $\hb-\hbo$ but keep the same notation).

Let us recall the $L^2$ condition. We denote by $\wt\jmath$ the inclusion $X^*\hto\wt X$. Then $\wcLd^0(H,\varh)$ is the subsheaf of $\wt\jmath_*L_\loc^1(H)$ consisting of sections which are holomorphic with respect to $\hbo$ and $L^2$ with respect to the metric $\varh$ on each compact set of the open set on which they are defined.

Given a local section $u$ of $\wt\jmath_*L_\loc^1(H)$ on $\wt X$, written as $\sum u_{\beta,\ell,k}(r,\theta) e^{\prime(\hbo)}_{\beta,\ell,k}$, it is a local section of $\wcLd^0(H,\varh)$ iff
\begin{equation}\label{eq:L2FH}
\big[(r,\theta)\ra u_{\beta,\ell,k}(r,\theta)\cdot r^{\ell_{\hbo}(q_{\beta,\imhb_o}+\beta)}\rL(r)^{\ell/2-1}e_\varh\big]\in L^2(d\theta\,dr/r),
\end{equation}
with $e_\varh=1,e^{\reel(\hbo/\ov{t'})},1+e^{\reel(\hbo/\ov{t'})}$ if $\varh=h,h_\hbo,h+h_\hbo$ (\cf \cite[p.~135]{errBibi01c} for the notation).

We define similarly $\wcLd^1(H,\varh)$ and $\wcLd^2(H,\varh)$ by asking moreover that $\omega_r,\omega_\theta$ have norm $\rL(r)$ and $\omega_r\wedge\omega_\theta$ has norm $\rL(r)^2$ (up to some constant depending on~$\hbo$). Therefore, a local section $v$ of $\wcLd^1(H,\varh)$ has coefficients $v^{(r)}_{\beta,\ell,k}$ and $v^{(\theta)}_{\beta,\ell,k}$ on $e^{\prime(\hbo)}_{\beta,\ell,k}\omega_r$ and $e^{\prime(\hbo)}_{\beta,\ell,k}\omega_\theta$ respectively, which satisfy \eqref{eq:L2FH} with $\rL(r)^{\ell/2}$ instead of $\rL(r)^{\ell/2-1}$. Similarly, a local section $w$ of $\wcLd^2(H,\varh)$ has coefficients $w_{\beta,\ell,k}$ on $e^{\prime(\hbo)}_{\beta,\ell,k}\omega_r\wedge\omega_\theta$ which satisfy \eqref{eq:L2FH} with $\rL(r)^{\ell/2+1}$ instead of $\rL(r)^{\ell/2-1}$.

\begin{lemme}[$\wt L^2$ Poincar\'e Lemma]\label{lem:Fpoin}
The complexes $\wcLd^{1+\cbbullet}(H,\varh,\wt\ccD_\hbo)$ ($\varh=h$, $h_\hbo$ or~$h+h_\hbo$) have cohomology in degree~$-1$ at most.
\end{lemme}

Keeping the notation of \cite[(5.3.7)]{errBibi01c}, the matrix $\wt\Theta'_\hbo$ of $\wt\ccD_\hbo$ in the basis $\bme^{\prime(\hbo)}$ can be decomposed as
\begin{gather*}
\wt\Theta'_\hbo=\wt\Theta'_{\hbo,\diag}+\Theta'_{\hbo,\nilp}+\Theta'_{\hbo,\pert},\\
\tag*{with}
\wt\Theta'_{\hbo,\diag}=\oplus_\beta\big[(q_{\beta,\imhb_o}+\beta)\star\hbo+1/t'\big]\id\,\dfrac{dt'}{t'}\\
\tag*{and}\Theta'_{\hbo,\nilp}=\big[\rY+P(0,\hbo)\big]\dfrac{dt'}{t'},\quad
\Theta'_{\hbo,\pert}=\big[P(t',\hbo)-P(0,\hbo)\big]\dfrac{dt'}{t'},
\end{gather*}
with $\rY=(\oplus_\beta \rY_\beta)$ (\cf \cite[Proof of Theorem 6.2.5]{errBibi01c}). We set $\rN_\hbo=\rY+P(0,\hbo)$.

Using \cite[Formula (6.2.7)]{errBibi01c}, we see as in \loccit that the $L^2$ condition on derivatives under~$\wt\ccD_\hbo$ can be replaced with the $L^2$ condition on derivatives under $\wt\ccD_{\hbo,\diag}$ (having matrix $\wt\Theta'_{\hbo,\diag}$): indeed, $\Theta'_{\hbo,\nilp}+\Theta'_{\hbo,\pert}$ sends $L^2$ sections to $L^2$ sections, when using the metric $\varh$.

Let $\theta_o\in S^1$, $r_1\in{}]0,r_0[$ and let $U={}]0,r_1[{}\times{}]\theta_o-\epsilon,\theta_o+\epsilon[$ be an open sector in $X^*$ with $\epsilon>0$ small enough so that $[\theta_o-\epsilon,\theta_o+\epsilon]$ contains at most one zero of $\cos(\theta+\arg\hbo)\cdot\sin(\theta+\arg\hbo)$ and this zero belongs to the interior of the interval. We denote by $\ov U$ its (compact) closure.

If $\cE^k$ denotes the sheaf of $C^\infty$ $k$-forms,
\[
\Gamma(\ov U,\cLd^k(H,\varh,\wt\ccD_\hbo))=L^2(\ov U,\cE^k_{\ov U}\otimes H,\varh,\wt\ccD_{\hbo,\diag})
\]
and the right-hand term is a Hilbert space, the norm being given by $\|\cdot\|_{2,\varh}+\|\wt\ccD_{\hbo,\diag}\cdot\nobreak\|_{2,\varh}$.

The proof will decompose in $3$ steps:
\begin{itemize}
\item
We first prove the lemma for the $L^2$ complex
\[
\big(L^2(\ov U,\cE^{1+\cbbullet}_{\ov U}\otimes H,\varh,\wt\ccD_{\hbo,\diag}),\wt\ccD_{\hbo,\diag}\big),
\]

\item
without changing the terms of the complex, we change the differential to $\wt\ccD_{\hbo,\diag}+\Theta'_{\hbo,\nilp}$ and prove the lemma by an extension argument,
\item
last, we change the differential to $\wt\ccD_{\hbo,\diag}+\Theta'_{\hbo,\nilp}+\Theta'_{\hbo,\pert}$, that we regard as a small perturbation of the previous one.
\end{itemize}

\begin{proof}[Proof of Lemma \ref{lem:Fpoin}, first step]
It is permissible to rescale the basis $\bme^{\prime(\hbo)}$, which therefore remains $L^2$-adapted (\cf \cite[\S6.2.b]{errBibi01c}), by multiplying each term $e^{\prime(\hbo)}_{\beta,\ell,k}$ by the function $e^{1/\hbo t'}t^{\prime-(q_{\beta,\imhb_o}+\beta)\star\hbo/\hbo}$ to obtain a basis $\wt{\bme^{\prime(\hbo)}}$, which is $\wt\ccD_{\hbo,\diag}$-flat. On the other hand, the $h$-norm of $e^{\prime(\hbo)}_{\beta,\ell,k}$ is equivalent, when $t'\to0$, to $\mtp^{\ell_\hbo(q_{\beta,\imhb_o}+\beta)}\Ltp^{\ell/2}$ up to a multiplicative constant depending on $\hbo$ (\cf \cite[Formula (5.3.6)]{errBibi01c}).

Therefore, the $h_\hbo$-norm of $\wt{e^{\prime(\hbo)}}_{\beta,\ell,k}$ is equivalent (up to a constant) to
\[
e^{\tfrac{1+\module{\hbo}^2}{\module{\hbo}r}\cos(\theta+\arg\hbo)}\cdot r^{-\tfrac{\beta''}{2}(\module{\hbo}+1/\module{\hbo})\sin\arg\hbo}\cdot\rL(r)^{\ell/2-1}.
\]
On the other hand, the $h$-norm is given by the same formula, where we replace $e^{\tfrac{1+\module{\hbo}^2}{\module{\hbo}r}\cos(\theta+\arg\hbo)}$ with $e^{\tfrac{1}{\module{\hbo}r}\cos(\theta+\arg\hbo)}$.

The proof of the vanishing of the higher cohomology sheaves in all three cases is then completely similar to that of \cite[Lemma~4.1]{errBibi98}.
\end{proof}

\begin{proof}[Proof of Lemma \ref{lem:Fpoin}, second step]
Consider the monodromy filtration of $\rN_\hbo$ and apply Step one to each graded piece. Use then an easy extension argument.
\end{proof}

\begin{proof}[Proof of Lemma \ref{lem:Fpoin}, third step]
We then apply to the complex of Hilbert spaces considered in Step two a standard perturbation argument, as the $L^2$-norm of $\Theta'_{\hbo,\pert}$ can be made small if $r_1$ is small (see \eg \cite[Lemma~2.68, p.~53]{errMochizuki07}).
\end{proof}

\subsubsection*{The complex $\wDR(\FccM_\hbo)$}

We extend the coefficients of $\FccM_\hbo$ to $\cA^\rmod_{\wt X}$ and we consider the corresponding de~Rham complex, that we denote by $\wDR(\FccM_\hbo)$. This is a complex on $\wt X$. Let us note that, as $\bR\rho_*\cA^\rmod_{\wt X}=\cO_X[t^{\prime-1}]$ and as $\FccM_\hbo$ is $\cO_X[t^{\prime-1}]$-flat (being locally free as such), we have $\bR\rho_*\wDR(\FccM_\hbo)=\DR\FccM_\hbo$.

\begin{lemme}[$\cA^\rmod_{\wt X}$-Poincar\'e lemma]\label{lem:APoinc}
The complex $\wDR(\FccM_\hbo)$ has cohomology in degree~$-1$ at most.
\end{lemme}

\begin{proof}
This is a particular case of a general result on irregular meromorphic connections, see \eg \cite[App.~1]{errMalgrange91}.
\end{proof}

\begin{lemme}[Comparison]\label{lem:comparisontilde}
The subsheaves $\cH^{-1}\wDR(\FccM_\hbo)$ and $\cH^{-1}\wcLd^{1+\cbbullet}(H,\varh,\wt\ccD_\hbo)$ ($\varh=h$, $h_\hbo$, or~$h+h_\hbo$) of $\wt\jmath_*j^{-1}\cH^{-1}\DR(\FccM_\hbo)$ coincide.
\end{lemme}

\begin{proof}
A $\wt\ccD_\hbo$-flat local section $u$ of $H$ takes the form $e^{1/\hbo t'}v$, where $v$ is a $\ccD_\hbo$-flat local section of $H$. Using for instance (5.3.6) and Remark 5.3.8(4) in \cite{errBibi01c}, one knows that the $h$-norm of $v$ grows exactly like $\mtp^b\Ltp^\nu$ for some $b\in\RR$ and some $\nu\in\frac12\ZZ$ when $t'\to0$. As the $h$-norm of $u$ is equal to $e^{\tfrac{1}{|\hbo|r}\cos(\theta+\arg\hbo)}\|v\|_h$, this norm is $L^2$ near $(\theta_o,\hbo)$ if and only if $\cos(\theta_o+\arg\hbo)<0$. The germ of $\cH^{-1}\wcLd^{1+\cbbullet}(H,h,\wt\ccD_\hbo)$ at $\theta_o$ is therefore~$0$ if $\cos(\theta_0+\arg\hbo)\geq0$, and consists of all flat local sections if $\cos(\theta_0+\arg\hbo)<0$.

Considering the metric $h_\hbo$ instead of $h$ will only replace $e^{\tfrac{1}{|\hbo|r}\cos(\theta+\arg\hbo)}$ with $e^{\tfrac{(1+|\hbo|^2)}{|\hbo|r}\cos(\theta+\arg\hbo)}$, so the argument is the same. The argument for $h+h_\hbo$ is also the same.

A similar argument shows that a $\wt\ccD_\hbo$-flat section has coefficients with moderate growth in the basis $\bme^{\prime(\hbo)}$ if and only if \hbox{$\cos(\theta_o+\arg\hbo)<0$} and, in such a case, any flat local section is a section of $\cH^{-1}\wDR(\FccM_\hbo)$.
\end{proof}

\begin{proof}[Proof of \eqref{enum:hhhbo}]
The assertion follows from Lemmas \ref{lem:Fpoin} and \ref{lem:comparisontilde} by taking $\bR\rho_*$. Let us note indeed that the complexes $\wcLd^{1+\cbbullet}$ are $c$-soft and that $\bR\rho_*\wcLd^{1+\cbbullet}=\cLd^{1+\cbbullet}$.
\end{proof}

\begin{proof}[Proof of \eqref{enum:DR2h}]
In order to prove \eqref{enum:DR2h}, we have to compare the complexes $\DR(\FccM_\hbo)$ and $\cLd^{1+\cbbullet}(H,h,\wt\ccD_\hbo)$. We will compare them with a third complex that we introduce now. We denote by $\Dbm_{\wt X}$ (\resp $\Dbm_X$) the sheaf on $\wt X$ (\resp $X$) of distributions on $X^*$ which can be lifted as distributions on $\wt X$ (\resp $X$). We have $\rho_*\Dbm_{\wt X}=\Dbm_X$. If $\Db_X$ is the sheaf of distributions on $X$, it is known that $\Dbm_X=\Db_X[t^{\prime-1}]$. We can define the complex on $\wt X$ of currents with moderate growth with values in $\FccM_\hbo$, that we denote by $\Dbmb_{\wt X}\otimes\rho^{-1}(\FccM_\hbo)$ and we have an inclusion $\wDR(\FccM_\hbo)\hto\Dbmb_{\wt X}\otimes\rho^{-1}(\FccM_\hbo)$. By an adaptation of the Dolbeault-Grothendieck theorem (\cf \cite[Prop.~II.1.1.7]{errBibi97}), the complex of moderate currents of type $(0,\cbbullet)$ with differential $d''$ is a resolution of $\cA_{\wt X}^\rmod$, hence the previous morphism is a quasi-isomorphism which becomes, after taking $\bR\rho_*$, the quasi-isomorphism $\DR(\FccM_\hbo)\to\Dbmb_X\otimes\FccM_\hbo$ (\cf \cite[\S2.c]{errBibi98}).

As the basis $\bme^{\prime(\hbo)}$ is $L^2$ adapted and as the $h$-norm of each element of this basis has moderate growth, we have a natural morphism from the $L^2$ complex to the complex of currents, that is, we have morphisms
\[
\wcLd^{1+\cbbullet}(H,h,\wt\ccD_\hbo)\hto\Dbmb_{\wt X}\otimes\rho^{-1}(\FccM_\hbo)\Hfrom{\sim}\wDR(\FccM_\hbo).
\]
From Lemma \ref{lem:comparisontilde} we conclude that the left morphism is a quasi-isomorphism, and finally, taking $\bR\rho_*$, we find quasi-isomorphisms
\[
\cLd^{1+\cbbullet}(H,h,\wt\ccD_\hbo)\Hto{\sim}\Dbmb_X\otimes\FccM_\hbo\Hfrom{\sim}\DR(\FccM_\hbo).
\]
Using now Lemma~\ref{lem:new3}, we find that the natural morphism
\begin{equation}\label{eq:22}
\DR(\FccM_{\hbo,\loc})_{(2)}\to\cLd^{1+\cbbullet}(H,h,\wt\ccD_\hbo)
\end{equation}
is a quasi-isomorphism.
\end{proof}

\section{Proof of the twistor property}\label{secapp:B}

In this section, it will be simpler to replace isometrically $(H,\Fou h,\Fou \ccD_\hbo))$, as defined in \S2.1 above, with $(H,h,\Lou\ccD_\hb)$, where $\Lou\ccD_\hb=e^{\ov t}\,\Fou\ccD_\hb e^{-\ov t}=\ccD_\hb-dt-\hb d\ov t$. We denote by $\Harm$ the space of harmonic sections in $\Gamma(\PP^1,\cLd^1(H,h,\Lou\ccD_\hbo))$. From the proof of Proposition~1 in the original text (as corrected above), we know that $\Harm$ does not depend on $\hbo$ when regarded as a subspace of $\Gamma(\PP^1,\cLd^1(H,h))$.

We denote by $\cP^1$ (\resp $\wt\cP^1$) the product $\PP^1\times\Omega_0$ (\resp $\wt\PP^1\times\Omega_0$), by $\rho$ the projection $\wt\cP^1\to\cP^1$ and by $\pr:\cP^1\to\Omega_0$ (\resp $\wt\pr=\rho\circ\pr:\wt\cP^1\to\Omega_0$) the natural projection. We define the $L^2$ sheaves on $\cP^1$ (\resp $\wt\cP^1$) in the same way as we did in \cite[\S6.2.b]{errBibi01c}. These sheaves are $\pr$-soft (\resp $\wt\pr$-soft) (\cf\cite[Def.~3.1.1]{errK-S90}). We thus have a natural morphism $\Harm\otimes_\CC\cO_{\Omega_0}\to \pr_*\cLd^1(\cH,h)$ constructed as in \cite[\S2.2.b]{errBibi01c}, and harmonic sections are in the kernel of $\Lou\ccD_\hb$ for any $\hb$, so the morphism takes values in $\pr_*\cLd^1(\cH,h,\Lou\ccD_\hb)$. Using the isometry given by the multiplication by $e^{-\hb\ov t}$, we find a natural morphism
\begin{equation}\label{eq:HarmL2hb}
\Harm\otimes_\CC\cO_{\Omega_0}\To{{}\cdot e^{-\hb\ov t}} \pr_*\cLd^1(\cH,h_\hb,\wt\ccD_\hb).
\end{equation}

We want to show that $\Harm$ is a lattice in $\bR^0\pr_*\DR\FcM$, and we will first find a morphism $\Harm\otimes_\CC\cO_{\Omega_0}\to\bR^0\pr_*\DR\FcM$.

\subsubsection*{The meromorphic $L^2$ de~Rham complex}
Let us first state an analogue of Lemma~3 in the original text. We denote by $\FcM_\loc$ the $\cR_{\cP^1}[*P]$-module obtained by localizing $\FcM$ at its singularities $P$. Note that, $\FcM_\loc$ coincides with $\FcM$ near $\infty$. The meromorphic $L^2$ de~Rham complex, with respect to the metric $h$, is denoted by $\DR(\FcM_\loc)_{(2),h}$. It is the sub-complex of $\DR\FcM_\loc$ defined by $L^2$ conditions with respect to $h$ for the sections and their derivatives. We have a natural morphism $\DR(\FcM_\loc)_{(2),h}\to\DR\FcM$: this is shown in \cite[\S6.2.a]{errBibi01c} at finite distance, and is clear near $\infty$.

\begin{lemme}\label{lem:DR2DR}
The natural morphism $\DR(\FcM_\loc)_{(2),h}\to\DR\FcM$ is a quasi-isomorphism.
\end{lemme}

\begin{proof}
This is \cite[Prop.~6.2.4]{errBibi01c} at finite distance and is proved as in Lemma~3 of the original text near $\infty$.
\end{proof}

\subsubsection*{The complex $\cF^{\protect\cbbullet}$}
As in the proof of Lemma \ref{lem:comparison}, we wish to work with moderate distributions near $\infty$, while keeping $L^2$ complexes at finite distance. We will denote by $X$ an open disc near $\infty$ in $\PP^1$ which contains no other singularity of $\FcM$ than $\infty$ and by $Y$ the complement of $\infty$ in $\PP^1$. Last, we set $Z=X\cap Y$. We will denote by $j_X:X\hto\PP^1$ the inclusion, and similarly for $j_Y$ and $j_Z$. We denote by the same letters the inclusion $\cX\hto\cP^1$, with $\cX=X\times\Omega_0$, etc.

We denote by $\Db_\cX$ the sheaf of distributions on $X$ and by $\Db_\cX^\an$ the sub-sheaf of distributions which are holomorphic with respect to $\hb$, \ie the kernel of $\ov\partial_\hb$. We denote by $(\Db_\cX^{\an,1+\cbbullet},\hb d'+d'')$ the sheaf of $\hb$-holomorphic currents on $\cX$ (we use the same rescaling on forms and currents as in \cite[\S0.3]{errBibi01c}). The Dolbeault-Grothendieck theorem implies that the complex of currents $(\Db_\cX^{\an,(k,0)},d'')$ is a resolution of $\Omega^k_\cX$. As $\FcM_{|\cX}$ is $\cO_\cX[*\infty]$-locally free (this follows from \cite[Lemma~5.4.1 and Lemma~3.4.1]{errBibi01c}) it is $\cO_\cX$-flat and $(\Db_\cX^{\an,(k,0)}\otimes_{\cO_\cX}\FcM_{|\cX},d'')$ is a resolution of $\Omega^k_\cX\otimes_{\cO_\cX}\FcM_{|\cX}$. Finally, we find that the natural morphism $\DR\FcM_{|\cX}\to\Db_\cX^{\an,1+\cbbullet}\otimes_{\cO_\cX}\FcM_{|\cX}$ is a quasi-isomorphism.

On the other hand, we have a morphism of complexes
\begin{equation}\label{eq:L2Db}
\cLd^{1+\cbbullet}(\cH,h,\wt\ccD_\hb)_{|\cX}\To{\iota}\Db_\cX^{\an,1+\cbbullet}\otimes_{\cO_\cX}\FcM_{|\cX}
\end{equation}
which, when restricted to $Z$, is a quasi-isomorphism. Indeed, on $Z$ this is clear. Near~$\infty$, this can bee seen by using the local $\cO_\cX[*\infty]$-basis $\bme^{\prime(\hbo)}$ of $\cM_\loc$ near $\infty$: this is a $L^2$-adapted basis and the $h$-norm of its elements has moderate growth near~$\infty$, locally uniformly with respect to $\hb$; this implies that a section of $\cLd(\cH,h)$ belongs to $\Db_\cX^\an\otimes_{\cO_\cX}\FcM_{|\cX}$. Let us check the compatibility of the differentials of the complexes. On $\cLd$, the derivative is not taken in the distributional sense on $\cX$, but only on $\cX^*=(X\moins\{\infty\})\times\Omega_0$. In other words, it is obtained by taking the derivative in the distributional sense on $\cX$ and then restricting to $\cX^*$. But the morphism $\iota$ is clearly compatible with this way of taking derivatives, as $\mtp$ acts in an invertible way on the right-hand side of \eqref{eq:L2Db}, hence any distribution supported on $\{\infty\}\times\Omega_0$ is annihilated by~$\iota$. (Let us notice that this point is exactly what prevents us from using distributions near singularities at finite distance, as $\FcM\neq\FcM_\loc$ near such a singular point.)

The complex $\cF^\cbbullet$ is defined by the exact sequence of complexes
\begin{multline*}
0\to j_{Z,!}\cLd^{1+\cbbullet}(\cH,h,\wt\ccD_\hb)_{|Z}\\
\To{(\id,-\iota)}j_{Y,!}\cLd^{1+\cbbullet}(\cH,h,\wt\ccD_\hb)_{|Y}\oplus j_{X,!}(\Db_\cX^{\an,1+\cbbullet}\otimes_{\cO_\cX}\FcM_{|\cX})\to\cF^\cbbullet\to0.
\end{multline*}
Let us note that each term in $\cF^\cbbullet$ is $\pr$-soft (\cf \cite[Prop.~2.5.7(ii) and Cor.~2.5.9]{errK-S90}).

\begin{lemme}\label{lem:DR2F}
We have a natural morphism of complexes $\DR(\FcM_\loc)_{(2),h}\to\cF^\cbbullet$ which is a quasi-isomorphism.
\end{lemme}

\begin{proof}
We use the exact sequence
\begin{multline*}
0\to j_{Z,!}j_Z^{-1}\DR(\FcM_\loc)_{(2),h}\\
\to j_{Y,!}j_Y^{-1}\DR(\FcM_\loc)_{(2),h}\oplus j_{X,!}j_X^{-1}\DR(\FcM_\loc)_{(2),h}\to \DR(\FcM_\loc)_{(2),h}\to0
\end{multline*}
to reduce the question to each of the open sets $X,Y,Z$. On $Y$, this is \cite[Th.~6.2.5]{errBibi01c}. On $Z$, this is easy, and on $X$, this follows from Lemma~\ref{lem:DR2DR}. The compatibility with the arrows in the previous exact sequences is easy.
\end{proof}

\begin{lemme}\label{lem:HarmF}
We have a natural morphism $\Harm\otimes_\CC\cO_{\Omega_0}\to\bR^0\pr_*\cF^\cbbullet=\cH^0(\pr_*\cF^\cbbullet)$.
\end{lemme}

\begin{proof}
Let us first note that the second equality comes from the $\pr$-softness of the terms in $\cF^\cbbullet$. Using \eqref{eq:L2Db}, we have a natural morphism $\cLd^{1+\cbbullet}(\cH,h,\wt\ccD_\hb)\to\cF^\cbbullet$. Therefore, it is enough to find a morphism
\begin{equation}\label{eq:HarmL2h}
\Harm\otimes_\CC\cO_{\Omega_0}\to\bR^0\pr_*\cLd^{1+\cbbullet}(\cH,h,\wt\ccD_\hb)=\cH^0\big(\pr_*\cLd^{1+\cbbullet}(\cH,h,\wt\ccD_\hb)\big).
\end{equation}
We have inclusions of $L^2$ complexes
\[
\cLd^{1+\cbbullet}(\cH,h,\wt\ccD_\hb)\Hfrom{\iota_h}\cLd^{1+\cbbullet}(\cH,h+h_\hb,\wt\ccD_\hb)\Hto{\iota_{h_\hb}}\cLd^{1+\cbbullet}(\cH,h_\hb,\wt\ccD_\hb).
\]

We will prove:
\begin{enumerate}\let\theenumi\theequation
\refstepcounter{equation}
\item\label{enum:nb0}
On some open neighbourhood $\nb(0)$ of $0$ in $\Omega_0$, the morphism \eqref{eq:HarmL2hb} factorizes through $\pr_*\iota_{h_\hb}$.
\refstepcounter{equation}
\item\label{enum:neq0}
On $\Omega_0\moins\{0\}$, the morphism $\iota_{h_\hb}$ is a quasi-isomorphism.
\end{enumerate}

This will be enough to conclude that we have a natural morphism
\[
\Harm\otimes_\CC\cO_{\Omega_0}\to\bR^0\pr_*\cLd^{1+\cbbullet}(\cH,h+h_\hb,\wt\ccD_\hb)=\cH^0\big(\pr_*\cLd^{1+\cbbullet}(\cH,h+h_\hb,\wt\ccD_\hb)\big),
\]
giving thus \eqref{eq:HarmL2h} by composing with $\bR^0\pr_*\iota_h$.
\end{proof}

\begin{proof}[Proof of \eqref{enum:nb0}]
By construction, $\Harm$ is a subspace of $\Gamma(\PP^1,\cLd^1(H,h,\Lou\ccD_\hb))$. We will use the following lemma, whose proof is due to the the referee of \cite{errBibi05} (note that S.~Szabo proves a similar result in \cite[Lemma 2.32]{errSzabo04}, with different methods however). If $f$ is a section of $H$ (\resp $\omega$ is a section of $H$ with values in $1$-forms), we will denote by $|f|_\varh$ (\resp $|\omega|_\varh$) the $\varh$-norm of $f$ (\resp the norm of $\omega$ with respect to $\varh$ and the norm induced by the Poincar\'e metric on $1$-forms, that we call the $\rP$-norm).

\begin{lemme}[Exponential decay of harmonic sections]\label{lem:decay}
For any $\omega\in\Harm$, there exists $C>0$ and a neighbourhood of~$\infty$ in~$X$ on which the $h$-norm of $\omega$ is bounded by $e^{-C|t|}$.
\end{lemme}

Once this lemma is proved, we obtain that $|e^{-\hb\ov t}\omega|_{h_\hb}=|\omega|_h\leq e^{-C|t|}$ for any $\omega\in\Harm$ on a suitable neighbourhood of $\infty$, hence $|\omega|_{h_\hb}\leq e^{-C|t|+\reel\hb\ov t}$. If $|\hb|$ is small enough, we thus get $|\omega|_{h_\hb}\leq e^{-C'|t|}$, and therefore $\omega$ is $L^2$ with respect to~$h_\hb$, as wanted.
\end{proof}

\begin{proof}[Proof of Lemma \ref{lem:decay}]
Let $\omega\in\Harm$. Then $\Lou\ccD_\hb\omega=0$ for any $\hb\in\Omega_0$, hence, if we set $\Lou\theta'_E=\theta'_E-dt$ and $\Lou\theta''_E=\theta''_E-d\ov t$, we have $(D''_E+\Lou\theta'_E)\omega=0$ and $(D'_E+\Lou\theta''_E)\omega=0$. We will now restrict the question near $\infty$ and we will work with the coordinate $t'$.

By the Dolbeault lemma for $\hbo=0$ (Lemma~4 of the original text corrected as in \S\ref{subsec:Dolbeault}), the complex $\cLd^{1+\cbbullet}(H,h,\wt\ccD_0)=\cLd^{1+\cbbullet}(H,h,(D''_E+\Lou\theta'_E))$ is quasi-isomorphic to $\DR\FccM_0$. Let us note that the germ of $\DR\FccM_0$ at $\infty$ is quasi-isomorphic to $0$, as $\Lou\theta'_E=t^{\prime-2}(\id+\cdots)$ is invertible on the germ $\FccM_0$ at $\infty$. Therefore, the germ \hbox{$\cLd^{1+\cbbullet}(H,h,(D''_E+\Lou\theta'_E))_\infty$} is quasi-isomorphic to $0$ and there exists a neighbourhood $X$ of $\infty$ and a section $f\in L^2(X,H,h)$ such that $(D''_E+\Lou\theta'_E)f=\omega$. Assume we prove $|f|_h\leq e^{-C'/\mtp}$ for some constant $C'>0$. Then, according to the moderate growth of $\Lou\theta'_E$, we will also have $|\Lou\theta'_Ef|_h\leq e^{-C''/\mtp}$ for some $C''>0$ on some neighbourhood of $\infty$, and thus the desired inequality for the $(1,0)$ part of~$\omega$. Arguing with a conjugate argument, we get the same kind of inequality for the $(0,1)$ part, hence the lemma.

Let us note that $(D'_E+\Lou\theta''_E)(D''_E+\Lou\theta'_E)f=(D'_E+\Lou\theta''_E)\omega=0$, hence $D'_ED''_Ef=-\Lou\theta''_E\Lou\theta'_Ef$. Since $D'_ED''_E+D''_ED'_E=-(\Lou\theta'_E\Lou\theta''_E+\Lou\theta''_E\Lou\theta'_E)$, we also get $D''_ED'_Ef=-\Lou\theta'_E\Lou\theta''_Ef$ (all these equalities are taken on $X^*$ in the distributional sense).

In particular, as $D'_ED''_E+\Lou\theta''_E\Lou\theta'_E$ is elliptic on $X^*$, $f$ is $C^\infty$ on $X^*$. If we set $\Lou\theta'_E=\Lou\Theta'_Edt'$ and $\Lou\theta''_E=\Lou\Theta''_Ed\ov{t'}$, $\Lou\Theta''_E$ is the $h$-adjoint of $\Lou\Theta'_E$. We then have on~$X^*$
\[
d'd''|f|_h^2=h(D'_ED''_Ef,\ov f)-h(D''_Ef,\ov{D''_Ef})+h(D'_Ef,\ov{D'_Ef})+h(f,\ov{D''_ED'_Ef}),
\]
so that, dividing by $dt'\wedge d\ov{t'}$ and using the previous relations, we find
\begin{equation}\label{eq:delta}
\partial_{t'}\partial_{\ov{t'}}|f|^2_h\geq|\Lou\Theta'_Ef|_h^2+|\Lou\Theta''_Ef|_h^2\geq C\mtp^{-4}|f|^2_h.
\end{equation}
This relation holds on $X^*$.

\begin{assertion}
The inequality \eqref{eq:delta} holds on $X$ in the weak sense, that is, for any nonnegative test function $\chi$ on $X$, and denoting by $d\vol_\rE$ the Euclidean volume $\itwopi dt'\wedge d\ov{t'}$,
\[
\int_X|f|^2_h(\partial_{t'}\partial_{\ov{t'}}\chi)\,d\vol_\rE\geq C\int_X\mtp^{-4}|f|^2_h\,\chi\,d\vol_\rE.
\]
\end{assertion}

\begin{proof}[Proof of the assertion]
Let us first note that $\mtp^{-2}|f|_h$ (hence also $|f|_h$) is in $L^2(d\vol_\rE)$, as $|\Lou\Theta'_Ef|_h$ is in $L^2(d\vol_\rP)$, where $d\vol_\rP=\mtp^{-2}\rL(t')^{-2}d\vol_\rE$ is the Poincar\'e volume, and $|\Lou\Theta'_Ef|_h\sim\mtp^{-2}|f|_h|dt'|_\rP$, with $|dt'|_\rP\sim \mtp\rL(t')$. In particular, $\|f\|_{h,\rP}<+\infty$. Similarly, if $\psi d\ov{t'}$ is the $(0,1)$ component of $\omega$, we have $\big|\int_Xh(D''_Ef,\ov{D''_Ef})\big|=2\pi\int_X|\psi|_h^2d\vol_\rE=2\pi\int_X|\psi|_h^2|d\ov{t'}|_\rP^2d\vol_\rP<\infty$, hence $|\psi|_h\in L^2(d\vol_\rE)$.

We now claim that $\big|\int_Xh(D'_Ef,\ov{D'_Ef})\big|<\infty$, that is, $\|D'_Ef\|_{h,\rP}<+\infty$. This follows from the acceptability (in the sense of \cite{errSimpson90}) of the Hermitian bundle $(H,D''_E,h)$ (as the Higgs field $\theta_E$ is tame). Indeed, the $\rP$-norm of the curvature $R(h)$ of $h$ is bounded near $\infty$. For any test function $\eta$ on $X^*$, we have (\cf \cite[(2.23)]{errMochizuki07})
\begin{align*}
\Big|\int_Xh(D'_E\eta,\ov{D'_E\eta})\Big|&\leq \Big|\int_Xh(D''_E\eta,\ov{D''_E\eta})\Big|+\Big|\int_X(\eta R(h),\eta d\vol_\rP)_{h,\rP}\,d\vol_\rP\Big|\\
&\leq\Big|\int_Xh(D''_E\eta,\ov{D''_E\eta})\Big|+\Big|\int_X|\eta|_h^2 |R(h)|_\rP\,|d\vol_\rP|_\rP \,d\vol_\rP\Big|\\
&\leq \Big|\int_Xh(D''_E\eta,\ov{D''_E\eta})\Big|+C\Big|\int_X|\eta|_h^2 d\vol_\rP\Big|,
\end{align*}
hence, $\|D'_E\eta\|_{h,\rP}\leq\|D''_E\eta\|_{h,\rP}+C\|\eta\|_{h,\rP}$. Since the Poincar\'e metric is complete near $\infty$, we can find a sequence of nonnegative test functions $\eta_n$ on $X^*$, which tend pointwise to $1$ in some punctured neighbourhood of $\infty$, such that $\eta_n\leq1$ and $|d\eta_n|_\rP\leq2^{-n}$ (see \eg \cite[Lemme~12.1]{errDemailly96}). Applying the previous result to $\eta_nf$, we find $\|\eta_nD'_Ef\|_{h,\rP}\leq\|\eta_nD''_Ef\|_{h,\rP}+(C+2^{-n+1})\|f\|_{h,\rP}$, hence the claim.

In order to end the proof of the assertion, it is enough to showing that the difference $\int_{\mtp\geq\epsilon}[|f|^2_h(\partial_{t'}\partial_{\ov{t'}}\chi)-(\partial_{t'}\partial_{\ov{t'}}|f|^2_h)\chi]\,d\vol_\rE$ tends to $0$ with $\epsilon$. It is then enough to find a sequence $\epsilon_n\to0$ such that $\int_{\mtp=\epsilon_n}|f|^2_h\,d\theta$, $\int_{\mtp=\epsilon_n}\partial_{t'}|f|^2_h\,d\theta$ and $\int_{\mtp=\epsilon_n}\partial_{\ov{t'}}|f|^2_h\,d\theta$ tend to $0$, and it is enough to checking that the integrals of $|f|^2_h,|\partial_{t'}|f|^2_h|,|\partial_{\ov{t'}}|f|^2_h|$ with respect to $d\theta\,dr/r$ is finite. For the first one, this follows from $\mtp^{-2}|f|_h\in L^2(d\vol_\rE)$. For the second one (and similarly the third one), we use that $|\varphi|_h$, $|\psi|_h$ and $\mtp^{-2}|f|_h$ belong to $L^2(d\vol_\rE)$.
\end{proof}

Once the assertion is proved, we can use the same trick (a variant of Ahlfors lemma) as in \cite{errSimpson90}. Let us remark first that, because $\partial_{t'}|f|^2_h$ and $\partial_{\ov{t'}}|f|^2_h$ are $L^1_\loc(d\vol_\rE)$ at $t'=0$, $|f|^2_h$ is continuous (and $C^\infty$ on $X^*$). Let us consider the auxiliary function $\exp(-C^{1/2}\mtp^{-1})$. A simple computation shows that $\partial_{t'}\partial_{\ov{t'}}\exp(-C^{1/2}\mtp^{-1})\leq C\mtp^{-4} \exp(-C^{1/2}\mtp^{-1})$. Let us then choose $\lambda>0$ such that $|f|^2_h\leq\lambda\exp(-C^{1/2}\mtp^{-1})$ in some neighbourhood of $\partial X$ and let $U\subset X$ be the open set where $|f|^2_h>\lambda\exp(-C^{1/2}\mtp^{-1})$. The previous inequalities show that $|f|^2_h-\lambda\exp(-C^{1/2}\mtp^{-1})$ is continuous and subharmonic in $U$. If $U$ is not empty then, at a boundary point of $U$ in $X$ we have $|f|^2_h=\lambda\exp(-C^{1/2}\mtp^{-1})$ and, by the maximum principle, we have $|f|^2_h-\lambda\exp(-C^{1/2}\mtp^{-1})\leq0$ on $U$, a contradiction.
\end{proof}

\begin{proof}[Proof of \eqref{enum:neq0}]
The proof is similar to that of Lemma~\ref{lem:Fpoin}. Let us work near $\hbo\in\Omega_0^*$. Using the $L^2$-adapted basis $\bme^{(\hbo)}$ we trivialize the bundle $\cH$ near $\hbo$. Given $\theta_o\in S^1$, we choose an open neighbourhood $\nb(\hbo)$ such that the choice of $r_1$ and $\epsilon$ in the proof of Lemma~\ref{lem:Fpoin} can be done uniformly with respect to $\hb\in\ov{\nb(\hbo)}$. Let $\rH(\nb(\hbo))$ denote the Banach space of continuous functions on $\ov{\nb(\hbo)}$ which are holomorphic in $\nb(\hbo)$. We then consider the complex whose terms are the $\bigoplus_{\beta,\ell,k}L^2\big(\ov U,\rH(\nb(\hbo)),\varh_{\beta,\ell,k},\wt\ccD_{\hb,\diag}\big)$ twisted by differential forms, where $\varh_{\beta,\ell,k}$ is $\|e^{(\hbo)}_{\beta,\ell,k}\|_{\varh,2}^2\varh$, and differential as in the three steps of the proof of Lemma~\ref{lem:Fpoin}.
We show as in Lemma~\ref{lem:Fpoin} that this complex has vanishing higher cohomology, and we obtain \eqref{enum:neq0}.
\end{proof}

\begin{proof}[Proof that $\Harm$ is a lattice]
From Lemmas \ref{lem:DR2DR}, \ref{lem:DR2F} and \ref{lem:HarmF} we get a morphism
\begin{equation}\label{eq:HarmDR}
\Harm\otimes_\CC\cO_{\Omega_0}\to \bR^0\pr_*\DR\FcM.
\end{equation}
As both terms are locally free $\cO_{\Omega_0}$-modules of the same rank, it will be an isomorphism as soon as its restriction to each fibre $\hb=
\hbo$ is an isomorphism of $\CC$-vector spaces. We will shorten the notation and denote by $_{|\hb=\hbo}$ the quotient by the image of $(\hb-\hbo)$.

For any complex $\cG^\cbbullet$ entering in the definition of the morphism \eqref{eq:HarmDR}, we have natural morphisms (with an obvious notation)
\[
(\bR^0\pr_*\cG^\cbbullet)_{|\hb=\hbo}\to\bR^0\pr_*(\cG^\cbbullet_{|\hb=\hbo})\to\bR^0\pr_*\cG^\cbbullet_\hbo.
\]

According to the exact sequence
\[
0\to\DR\FcM\To{\hb-\hbo}\DR\FcM\to\DR\FccM_\hbo\to0,
\]
and since each of these complexes have hypercohomology in degree $0$ at most, the natural morphism $(\bR^0\pr_*\DR\FcM)_{|\hb=\hbo}\to\bH^0(\PP^1,\DR\FccM_\hbo)$ is an isomorphism.

As a consequence, it is enough to prove that, for any $\hbo\in\Omega_0$, the morphism $\Harm\to\bH^0(\PP^1,\DR\FccM_\hbo)$ constructed as \eqref{eq:HarmDR} by fixing $\hb=\hbo$, is an isomorphism. Let us recall how it is constructed, by considering the following commutative diagram:
\[
\xymatrix@=5mm{
\Harm\ar[r]^-\sim&H^0(\PP^1,\cLd^{1+\cbbullet}(H,h_\hbo,\wt\ccD_\hbo)\\
&H^0(\PP^1,\cLd^{1+\cbbullet}(H,h+h_\hbo,\wt\ccD_\hbo)\ar[u]_-\wr \ar[d]^-\wr\\
&H^0(\PP^1,\cLd^{1+\cbbullet}(H,h,\wt\ccD_\hbo)\ar[d]_-c&\ar[l]_-a \bH^0(\PP^1,\DR(\FccM_{\hbo,\loc})_{(2)})\ar[d]^-\wr\\
&H^0(\PP^1,\cF^\cbbullet_\hbo)&\bH^0(\PP^1,\DR\FccM_\hbo)\ar[l]^-\sim_-b
}
\]
Then \eqref{eq:HarmDR}$_\hbo$ is obtained by factorizing through $H^0(\PP^1,\cF^\cbbullet_\hbo)$ and $b^{-1}$. On the other hand, we know that $a$ is an isomorphism (this is \eqref{eq:22} if $\hbo\neq0$ and Lemma~4 of the original text as corrected in \S\ref{subsec:Dolbeault} if $\hbo=0$). Therefore, $c$ is also an isomorphism.
\end{proof}

\begin{proof}[End of the proof of the twistor property]
The proof is done as in \cite[p.~53]{errBibi01c}, where we use the $L^2$ complex instead of the $C^\infty$ de~Rham complex.
\end{proof}

\providecommand{\bysame}{\leavevmode\hbox to3em{\hrulefill}\thinspace}
\providecommand{\MR}{\relax\ifhmode\unskip\space\fi MR }
\providecommand{\MRhref}[2]{%
\href{http://www.ams.org/mathscinet-getitem?mr=#1}{#2}
}

\providecommand{\href}[2]{#2}

\makeatletter
\let\enddoc@text
\oldenddoc@text
\makeatother

\begin{thebibliography}{10}

\bibitem{Biquard97}
{\scshape O.~Biquard} -- {\og {Fibr\'es de {H}iggs et connexions int\'egrables:
  le cas logarithmique (diviseur lisse)}\fg}, \emph{Ann. scient. {\'E}c. Norm.
  Sup. {$4^{\rm e}$} s{\'e}rie} \textbf{30} (1997), p.~41--96.

\bibitem{Borelal87}
{\scshape A.~Borel} (\smfedname) -- \emph{Algebraic {$\mathcal{D}$}-modules},
  Perspectives in Math., vol.~2, Boston, Academic Press, 1987.

\bibitem{CIMPA90-1}
{\scshape {\relax Ph}.~Maisonobe {\normalfont \smfandname} C.~Sabbah}
  (\smfedsname) -- \emph{{$\mathcal{D}$-modules coh\'erents et holonomes}}, Les
  cours du CIMPA, Travaux en cours, vol.~45, Hermann, Paris, 1993.

\bibitem{Malgrange91}
{\scshape B.~Malgrange} -- \emph{{\'E}quations diff\'erentielles {\`a}
  coefficients polynomiaux}, Progress in Math., vol.~96, Birkh{\"a}user, Basel,
  Boston, 1991.

\bibitem{Bibi90}
{\scshape C.~Sabbah} -- {\og Introduction to algebraic theory of linear systems
  of differential equations\fg}, in \emph{{\'E}l\'ements de la th\'eorie des
  syst\`emes diff\'erentiels} \cite{CIMPA90-1}, p.~1--80.

\bibitem{Bibi98}
\bysame , {\og {Harmonic metrics and connections with irregular
  singularities}\fg}, \emph{Ann. Inst. Fourier (Grenoble)} \textbf{49} (1999),
  p.~1265--1291.

\bibitem{Bibi00}
\bysame , \emph{{D\'eformations isomonodromiques et vari\'et\'es de
  Frobenius}}, Savoirs Actuels, CNRS~{\'E}ditions \& EDP~Sciences, Paris, 2002.

\bibitem{Bibi01}
\bysame , {\og {Polarizable twistor $\mathcal{D}$-modules}\fg},
  pr\'epublication, 218 pages, disponible \`a
  \url{http://math.polytechnique.fr/cmat/sabbah/articles.html}, 2004.

\bibitem{Simpson90}
{\scshape C.~Simpson} -- {\og Harmonic bundles on noncompact curves\fg},
  \emph{J.~Amer. Math. Soc.} \textbf{3} (1990), p.~713--770.

\bibitem{Simpson92}
\bysame , {\og Higgs bundles and local systems\fg}, \emph{Publ. Math. Inst.
  Hautes {\'E}tudes Sci.} \textbf{75} (1992), p.~5--95.

\bibitem{Simpson97}
\bysame , {\og Mixed twistor structures\fg}, Pr\'epublication Universit\'e de
  Toulouse \& \url{arXiv:math.AG/9705006}, 1997.

\bibitem{Szabo04}
{\scshape S.~Szabo} -- {\og {Nahm transform of meromorphic integrable
  connections on the Riemann sphere}\fg}, in preparation, 2004.

\bibitem{Zucker79}
{\scshape S.~Zucker} -- {\og {Hodge theory with degenerating coefficients:
  {$L_2$}-cohomology in the Poincar\'e metric}\fg}, \emph{Ann. of Math.}
  \textbf{109} (1979), p.~415--476.

\end{thebibliography}

\begin{thebibliography}{10}

\bibitem{errDemailly96} J.-P. Demailly, \emph{{Th\'eorie de Hodge $L^2$ et th\'eor\`emes d'annulation}}, Introduction {\`a} la th\'eorie de Hodge, Panoramas \& Synth\`eses, vol.~3, Soci{\'e}t{\'e} Math{\'e}matique de France, 1996, pp.~3--111.

\bibitem{errK-S90} M.~Kashiwara and P.~Schapira, \emph{{Sheaves on Manifolds}}, Grundlehren der mathematischen Wissenschaften, vol. 292, Springer-Verlag, 1990.

\bibitem{errMalgrange91} B.~Malgrange, \emph{{\'E}quations diff\'erentielles {\`a} coefficients polynomiaux}, Progress in Math., vol.~96, Birkh{\"a}user, Basel, Boston, 1991.

\bibitem{errMochizuki07} T.~Mochizuki, \emph{{Asymptotic behaviour of tame harmonic bundles and an application to pure twistor $D$-modules}}, vol. 185, Mem. Amer. Math. Soc., no. 869-870, American Mathematical Society, Providence, RI, 2007.

\bibitem{errBibi98} C.~Sabbah, \emph{{Harmonic metrics and connections with irregular singularities}}, Ann. Inst. Fourier (Grenoble) \textbf{49} (1999), 1265--1291.

\bibitem{errBibi97} \bysame, \emph{{{\'E}quations diff\'erentielles {\`a} points singuliers irr\'eguliers et ph\'enom\`ene de Stokes en dimension {$2$}}}, Ast{\'e}risque, vol. 263, Soci{\'e}t{\'e} Math{\'e}matique de France, Paris, 2000.

\bibitem{errBibi05} \bysame, \emph{{Fourier-Laplace transform of a variation of polarized complex Hodge structure}}, arXiv: \url{math.AG/0508551}, 36 pages, 2005.

\bibitem{errBibi01c} \bysame, \emph{{Polarizable twistor $\mathcal{D}$-modules}}, Ast{\'e}risque, vol. 300, Soci{\'e}t{\'e} Math{\'e}matique de France, Paris, 2005.

\bibitem{errSimpson90} C.~Simpson, \emph{Harmonic bundles on noncompact curves}, J.~Amer. Math. Soc. \textbf{3} (1990), 713--770.

\bibitem{errSzabo04} S.~Szabo, \emph{{Nahm transform of meromorphic integrable connections on the Riemann sphere}}, Ph.D. thesis, Universit\'e Louis Pasteur, Strasbourg, juillet 2005, arXiv: \url{0511471v1}.

\end{thebibliography}
\end{document}